\documentclass[11pt]{article}
        \usepackage{amssymb,amsmath,bm}
        \usepackage[colorlinks=true, urlcolor=blue,linkcolor=blue, citecolor=blue]{hyperref}
        \usepackage{mathrsfs,verbatim}
        \usepackage{caption,float}
        \usepackage{constants,color}
        \usepackage{enumerate}
        \usepackage{bbm}
        \usepackage{appendix}
        \usepackage{graphicx, caption,subcaption}
        \usepackage{graphicx}  
        \usepackage{booktabs}
        \usepackage{tikz}
        \usetikzlibrary{shapes.geometric, calc}
\begin{document}
        \newcommand{\bea}{\begin{eqnarray}}
        \newcommand{\ena}{\end{eqnarray}}
        \newcommand{\beas}{\begin{eqnarray*}}
        \newcommand{\enas}{\end{eqnarray*}}
        \newcommand{\beq}{\begin{equation}}
        \newcommand{\enq}{\end{equation}}
        \def\qed{\hfill \mbox{\rule{0.5em}{0.5em}}}
        \newcommand{\bbox}{\hfill $\Box$}
        \newcommand{\ignore}[1]{}
        \newcommand{\ignorex}[1]{#1}
        \newcommand{\wtilde}[1]{\widetilde{#1}}
        \newcommand{\mq}[1]{\mbox{#1}\quad}
        \newcommand{\bs}[1]{\boldsymbol{#1}}
        \newcommand{\qmq}[1]{\quad\mbox{#1}\quad}
        \newcommand{\qm}[1]{\quad\mbox{#1}}
        \newcommand{\nn}{\nonumber}
        \newcommand{\Bvert}{\left\vert\vphantom{\frac{1}{1}}\right.}
        \newcommand{\To}{\rightarrow}
        \newcommand{\supp}{\mbox{supp}}
        \newcommand{\law}{{\cal L}}
        \newcommand{\Z}{\mathbb{Z}}
        \newcommand{\di}{\mathrm{DI}}
                \newcommand{\bin}{\mathrm{Bin}}
        \newcommand{\ci}{\mathrm{CI}}
        \newcommand{\var}{\mathrm{Var}}
        \newcommand{\cov}{\mathrm{Cov}}

        \newcommand{\ucolor}[1]{\textcolor{blue}{#1}}  
        \newcommand{\ucomm}[1]{\marginpar{\tiny\ucolor{#1}}}  
        
        \newcommand{\bcolor}[1]{\textcolor{red}{#1}}  
        \newcommand{\bcomm}[1]{\marginpar{\tiny\bcolor{#1}}}  

        \newtheorem{theorem}{Theorem}[section]
        \newtheorem{corollary}{Corollary}[section]
        \newtheorem{conjecture}{Conjecture}[section]
        \newtheorem{proposition}{Proposition}[section]
        \newtheorem{lemma}{Lemma}[section]
        \newtheorem{definition}{Definition}[section]
        \newtheorem{example}{Example}[section]
        \newtheorem{remark}{Remark}[section]
        \newtheorem{case}{Case}[section]
        \newtheorem{condition}{Condition}[section]
        \newcommand{\pf}{\noindent {\it Proof:} }
        \newcommand{\proof}{\noindent {\it Proof:} }
        \frenchspacing
        

        \title{\bf Analysis of Clustering and Degree Index in Random Graphs and Complex Networks} 
        \author{\"{U}m\.{i}t I\c{s}lak\footnote{Bo\u{g}azi\c{c}i University, Mathematics Department, Istanbul, Turkey. email: umit.islak1@bogazici.edu.tr} \hspace{0.05in}\footnote{Middle East Technical University, Institute of Applied Mathematics, Ankara, Turkey. email: uislak@metu.edu.tr}  \quad Barış Yeşiloğlu\footnote{Bo\u{g}azi\c{c}i University, Mathematics Department, Istanbul, Turkey. email: baris.yesiloglu@bogazici.edu.tr} }\vspace{0.25in}
        \maketitle
        
         \begin{abstract}
        The purpose of this paper is to analyze the degree index and the clustering index in random graphs. The degree index in our setup is a certain measure of degree irregularity whose basic properties are well studied in the literature, and the corresponding theoretical analysis in a random graph setup turns out to be tractable. On the other hand, the clustering index, based on a similar reasoning, is first introduced in this manuscript. Computing exact expressions for the expected clustering index turns out to be more challenging even in the case of Erd\H{o}s-R\'enyi graphs, and our results are on obtaining relevant upper bounds. These are also complemented with observations based on Monte Carlo simulations. Besides the Erd\H{o}s-R\'enyi case, we also do simulation-based analysis for random regular graphs, the Barab\'asi-Albert model and the  Watts-Strogatz model. 
        

        \textbf{MSC2020 subject classifications:} 05C80, 05C82, 05C07, 60C05. 
        \end{abstract}
        

        \section{Introduction}

        Since the introduction of Erd\H{o}s-R\'enyi graphs around 1960 (\cite{er:59}, \cite{er:60}), random graphs have gained substantial importance in both applied and theoretical fields. Meanwhile, various statistics of interest for the underlying networks (or graphs) have been studied, including the average degree, clustering coefficient, efficiency, and modularity, among several others. See, for example, \cite{chenetal:22}, \cite{hofstadt:24}, \cite{kl:13} and \cite{siewetal:16} for relevant definitions and certain applications. Due to the development of graph-based  machine learning algorithms in recent years, these statistics have become even more popular, especially in classification problems \cite{sz:16}. 
        The purpose of this paper is to study two other types of graph characteristics in random graphs.

        In order to describe the two characteristics that we study, 
        let $\mathcal{G} = (V,E)$ be a graph where $V = \{1,2,\ldots,n\}$ for some $n \in \mathbb{N}$. Let $N(k)$ be the set of neighbors of $k \in V$ in $\mathcal{G}$. We denote the degree of $k$ by $d_k$. With this notation, we define the  $\alpha$-\textit{level} \textit{degree index} of $\mathcal{G}$ to be $$\di_{\alpha} (\mathcal{G} )= \sum_{1 \leq i < j \leq n} |d_i - d_j|^{\alpha}.$$ Here, $\alpha > 0$ is a real number, which is chosen to be in $\{1,2\}$ in this manuscript. Also, let us note that $\sum_{i < j}$ will mean $\sum_{1 \leq i < j  \leq n}$ below. Noting that $\di_{\alpha}(\mathcal{G}) = 0$ means that the graph under study is regular, this quantity is considered as  a measure of irregularity of the graph \cite{west}. On the other hand, in networks with high disparity in node connectivity, the index $\di_{\alpha}(\mathcal{G})$ would yield large values. The index overall provides a global measure of degree variability in the network and may, for example, help differentiate between an Erd\H{o}s-R\'enyi graph and a Barab\'asi-Albert network even if they have the same number of nodes and edges. 
        
        A definition of degree irregularity that has found significant use in the literature is the \textit{Albertson index} introduced in \cite{albertson:11}, in which the corresponding definition is similar to ours, but the sum under consideration is over the edge set instead of the vertices. Although the irregularity used in the present manuscript is quite natural to consider, the special case  $\di_{1}(\mathcal{G})$  was relatively recently introduced in \cite{abdo:14}, where the authors call it the \textit{total irregularity} of the graph. In addition, a form of the special case of $\di_{2}(\mathcal{G})$, in which the corresponding sum is over the edge set instead of the vertex set, was proposed in \cite{abdo:18}. Aside from these, there are other notions of graph degree irregularity; see \cite{abdo:19} for a partial list and \cite{ali:21} for a broad look at the topic. 
        
        The second graph characteristic of interest for us will be the clustering index, which, as far as we know, is not studied previously in the literature. Before introducing this index, recall that for given $i \in V$, when  $d_i \geq 2$,  the \textit{local clustering coefficient} of $i $ is defined to be $$C(i) = \frac{1}{ d_i (d_i - 1)} \sum_{j, k \in N(i): j  \neq k} \mathbf{1} (j \leftrightarrow k).$$ In the case where $d_i \in \{0,1\}$, we set $C(i) = 0$. In general, the clustering coefficient of a node $i \in V(\mathcal{G})$ is a measure of the likelihood of the neighbors of $i$ to be neighbors among themselves \cite{ws:98}. The clustering coefficient is extensively studied in the literature, and several variations and generalizations are proposed; see \cite{kb:19}, \cite{sv:05}  and \cite{Haoetal:18} for some exemplary work. The clustering coefficient is also studied in a random graph setup. For instance, \cite{lsy:17} analyzes the clustering coefficient in Erd\H{o}s-R\'enyi graphs and certain random regular graphs. \cite{ghz:13} does an investigation in the case of a generalized small world model. On the other hand, \cite{Haoetal:18} discusses a natural generalization of the standard clustering coefficient, and analyzes in both Erd\H{o}s-R\'enyi and small world settings. 

        Now, similar to the degree index,  for $\alpha >0$, we define  the $\alpha$-\textit{level clustering index} of $\mathcal{G}$ to be 
        $$\ci_{\alpha} (\mathcal{G}) = \sum_{i < j} |C(i) - C(j)|^{\alpha}. $$
        Here again $\alpha$ in our case is either 1 or 2. As noted above, although such an index seems natural to define, it is first introduced in this manuscript to the best of our knowledge. Before continuing with our contributions, let us now include some discussion of this newly introduced index. 

        The clustering index $\ci_{\alpha} (\mathcal{G})$ captures the heterogeneity in local clustering across the graph, offering insights that traditional global or average measures may miss. For instance, it can reveal whether certain nodes are embedded in tightly-knit communities where others are more isolated. This distinction is crucial in real-world networks, where some nodes may form dense subgraphs, while others may lie on the periphery with sparse neighborhoods \cite{bm:00}. In many systems, the core nodes tend to have high clustering, whereas the peripheral nodes exhibit low clustering. A high clustering index reflects this inequality, which could remain hidden when only considering global or average clustering coefficients.

Moreover, in networks with community structure (\cite{r:04}), nodes within the same community often display similar clustering levels. Hence, a low clustering index may indicate structural homogeneity, while a high value can signal the presence of multiple structurally distinct communities. As a last note, when analyzing temporal networks, the global clustering might remain stable over time, even as the roles of nodes shift. In contrast, the clustering index is sensitive to such changes in clustering variability, making it a possible tool for detecting structural transitions or evolving roles in dynamic graphs.

        Moving to the contributions of this article, our primary goal is to analyze $\di_{\alpha} (G)$ and $\ci_{\alpha} (G)$ when $G$ is an Erd\H{o}s-R\'enyi graph. For such random graphs, the degree index case turns out to be relatively straightforward, and precise formulas for the expected degree index are obtainable. In particular, the results in Section \ref{sec:degindex} tell that when $G$ is an Erd\H{o}s-R\'enyi graph with $n$ nodes and a fixed attachment probability $p$, we have, $$\mathbb{E}[\di_1(\mathcal{G})]  \sim \frac{2}{\sqrt{\pi}} \binom{n}{2} \sqrt{(n-2) p (1 - p)} \quad \text{and} \quad \mathbb{E}[\di_2(\mathcal{G})] = 6\binom{n}{3}p(1-p).$$  
        
        On the other hand, the computations for the case of the expected clustering index are more involved, and we are not able to obtain exact formulas. However, for the same random graph family, we show that the bounds $$\mathbb{E}[\ci_1(\mathcal{G})] \leq K_1 n \quad   \text{and} \quad \mathbb{E}[\ci_2(\mathcal{G})] \leq K_2,$$ hold for every $n$,  for some constants $K_1, K_2$ depending on $p$, but not on $n$. Our heuristic arguments yield that matching lower bounds are true for both cases, but we have not verified these rigorously yet. Moreover, we will complement our theoretical results by making observations on the two indices of interest for other random graph models via Monte Carlo simulations. The models we study include random regular graphs, the Barab\'asi-Albert model, and the Watts-Strogatz model, and relevant simulations and discussions will be in Section \ref{sec:otherrgm}. 

    The remainder of the paper is organized as follows. The next section contains some basic observations on the clustering index, which can be easily adapted to the degree index. This section also includes a brief discussion on the comparison of the extremal cases of the degree index and the clustering index. Afterwards, in Section \ref{sec:clusind} we analyze the clustering index in Erd\H{o}s-R\'enyi graphs, and provide upper bounds for the expected clustering index for both $\alpha = 1$ and $2$. Later, in Section \ref{sec:degindex}, a similar study is done for the degree index. Section \ref{sec:otherrandom} is then devoted to a simulation-based study for the indices under study for random regular graphs, the Barab\'asi-Albert model, and the Watts-Strogatz model. Lastly, we conclude the paper in Section \ref{sec:conclusion} with some discussions on possible future work. 
         
        \section{Some basic observations on the clustering index}\label{sec:basicobs}
        
        In this section, we make some elementary observations on the clustering index, beginning with the analysis of the extremal values. Below, whenever it is clear from the context, we assume that the node set of the underlying graph is $V = \{1,2,\ldots,n\}$ for some $n \in \mathbb{N}$. The following elementary lemma will be used to determine the extremal values of $\ci_{\alpha} (\mathcal{G})$.  

        \begin{lemma}\label{prop:elt}
            Let $a_1, a_2, \ldots, a_n$ be real numbers on $[0,1]$. Then for any $\alpha \geq 1$, we have,
        $$
        \sum_{1 \leq i<j \leq n}\left|a_i-a_j\right|^{\alpha} \leq \frac{n^2}{4}.
        $$
        \end{lemma} 
        
        \textbf{Proof.} Clearly, it suffices to show that $\sum_{1 \leq i<j \leq n}\left|a_i-a_j\right| \leq n^2 /4$. Now, without loss of generality, we assume $0 \leq a_1 \leq \cdots \leq a_n \leq 1$, and define $b_i:=a_i-a_{i-1}$ for every $i \in\{2, \cdots, n\}$.
        Then, for every $j>i$, we have $\left|a_i-a_j\right|=a_j-a_i=\sum_{k=i+1}^j b_k$. Therefore,
        $$
        S:=\sum_{1 \leq i<j \leq n}\left|a_i-a_j\right|=\sum_{1 \leq i<j \leq n} \sum_{k=i+1}^j b_k=\sum_{k=2}^n(k-1)(n+1-k) b_k.
        $$
        The reason for the last equality is that for the terms containing $b_k$, we must have $i \leq k-1$ and $j \geq k$. The number of such pairs is $(k-1)(n+1-k)$. The maximum value $(k-1)(n+1-k)$ can take is $\frac{n^2}{4}$. Therefore, $S \leq \frac{n^2}{4} \sum_{k=2}^n b_k=$ $\frac{n^2}{4}\left(a_n-a_1\right) \leq \frac{n^2}{4}$. \hfill $\square$

Now, the next corollary provides the extremal values of $\ci_{\alpha} (\mathcal{G})$. The lower bound in the corollary is clear, and the upper bound follows from Lemma \ref{prop:elt}. 
          
        \begin{corollary}\label{prop:ciub}
        For any graph $\mathcal{G} = (V,E)$  and for any $\alpha \geq 1$, we have $$0 \leq \ci_{\alpha} (\mathcal{G}) \leq \frac{n^2}{4}.$$
        \end{corollary}

        \begin{example}\label{eg:1} In this example, we discuss some elementary graphs for which the bounds stated in Corollary \ref{prop:ciub} are attained.
        Let $\alpha \geq 1$ be arbitrary. The lower bound in Corollary \ref{prop:ciub} is attained, for example, for the following graphs:
        \begin{itemize}
            \item The null graph and the complete graph;
            \item Any forest, and in particular, any tree.
        \end{itemize}

      As a simple example where the clustering index is maximized, let $m \geq 2$, $n =2m$, and consider the graph that is the graph union of an $m$-null set, and $K_m$ (complete graph with $m$ vertices). Figure \ref{fig:eg21} exemplifies such a construction when $n=8$.

      Here, the only contribution to the sum defining $\ci_{\alpha}$ is from the cases where a node from $K_m$ and another node from the $m$-null set are considered. Noting that the clustering coefficient of the node from the $m$-null set is 0 by definition, the observation that $\ci_{\alpha}$ takes the value $n^2/4$ follows immediately.
\begin{figure}[H]
    \centering
    \begin{tikzpicture}
  \begin{scope}[xshift=0cm, yshift=1cm]
    \foreach \i in {1,2,3,4} {
      \coordinate (v\i) at (45+360/4*\i:1.5);
      \fill[red] (v\i) circle(2pt);
    }

    \foreach \i in {1,2,3,4} {
      \foreach \j in {\i,...,4} {
        \ifnum\i<\j
          \draw[blue, thick] (v\i) -- (v\j);
        \fi
      }
    }

    \node at (0, -2) {Complete Graph $K_4$};
  \end{scope}

  \begin{scope}[xshift=5cm, yshift=1cm]
    \foreach \i in {1,2,3,4} {
      \coordinate (u\i) at (45+360/4*\i:1.5);
      \fill[green!60!black] (u\i) circle(2pt);
    }

    \node at (0, -2) {Null Graph $N_4$};
  \end{scope}
    \end{tikzpicture}
    \caption{A graph with $8$ vertices where $\ci_{\alpha} = 8^2/4 = 16$.}
    \label{fig:eg21}
\end{figure}
      
  \noindent    Another related example will be in Example \ref{eg:poly} where extremal situations in clustering and degree indices are compared and discussed.
         \end{example}

        Noting that the observations in Corollary \ref{prop:ciub} can be easily adapted to the degree index, we continue with the following result which provides a way to express the clustering index in terms of the order statistics of the corresponding local clustering coefficients.   
        
          \begin{lemma}\label{prop:ciord} (Expression in terms of order statistics)
        Let $\mathcal{G} = (V,E)$ be some graph.  Order the local clustering coefficients as $C_{(1)} \leq C_{(2)} \leq \cdots \leq C_{(n)}$. Then, $$\sum_{i < j} |C(i) - C(j)| = \sum_{i = 2}^{n - 1} (2 i - 1 - n) C_{(i)} + (n  - 1) (C_{(n)} - C_{(n-1)}).$$ 
          \end{lemma}
        
        \textbf{Proof.} Following the proof of  Lemma \ref{prop:elt} with $a_i = C_{(i)} $ gives 
        \begin{eqnarray*}
            \sum_{i < j} |C(i) - C(j)| = \sum_{k=2}^n (k - 1)(n+1- k) (C_{(k)}- C_{(k-1)}).
        \end{eqnarray*}
        Expanding the right hand side as $\sum_{k=2}^n (k - 1)(n+1- k) C_{(k)} - \sum_{k=2}^n (k - 1)(n+1- k) C_{(k-1)}$, changing the index in the second sum properly, and doing some elementary manipulations yield the result. \hfill $\square$

         Let us also examine extremal cases in random graphs. Clearly, any random tree model provides a random graph model for which the lower bound in Corollary \ref{prop:ciub} is achieved. The following example, which is a randomized version of the construction given in Example \ref{eg:1}, provides a simple random graph model where the upper bound in Corollary \ref{prop:ciub} is achieved asymptotically.

        \begin{example}
        Begin with two empty sets of nodes $V_1$ and $V_2$ at time 0. At each time $t \in \mathbb{N}$, a newly generated node $v_t$ is inserted in $V_1$ or $V_2$ with probabilities $p$ and $1- p$, respectively. If $v_t$ is inserted in $V_1$, we do not attach it to any present nodes. If it is inserted in $V_2$, then we attach it to all vertices present in $V_2$. Then, at time $n \in \mathbb{N}$, the vertices in $V_1$ form a null graph, and the vertices in $V_2$ form a complete graph.
        
        Letting $N$ be the number of vertices in $V_1$ at time $n$, $N$ is binomially distributed with parameters $n$ and $p$, and we  have,
        \begin{eqnarray*}
        \mathbb{E} \left[ \sum_{i < j} |C(i) - C(j)|   \right]&=& \mathbb{E}\left[ \mathbb{E} \left[ \sum_{i < j} |C(i) - C(j)| \, \bigg| \, N \right] \right] = \mathbb{E}[N (n - N)] = n \mathbb{E}[N] - \mathbb{E}[N^2] \\ &=& n^2 p - n p (1- p) - n^2 p^2  =  (n^2 -n) p (1 -p). 
        \end{eqnarray*}
        When $p = 1/2 $, we have $\mathbb{E} [\ci_1 (\mathcal{G}) ]= \frac{1}{4} (n^2 -n).$  In particular, for large $n$,  $\mathbb{E} [\ci_1 (\mathcal{G}) ] \approx \frac{n^2}{4}$.  
        \end{example}

        In the following last example, we provide a comparison of the clustering index and the degree index. In particular, we are interested in cases where the degree index is high and the clustering coefficient is low, and vice versa. 

        \begin{example}\label{eg:poly}
        Consider a collection of disjoint polygons, i.e. a graph where each node is contained in a unique polygon. $C(i) = 1$ if $i$ is contained in a triangle, and $C(i) = 0$ if $i$ is contained in a polygon other than a triangle. The degree of any node in this graph is 2. Assuming we have $k$ vertices contained in some triangle and $n-k$ vertices contained in a polygon which is not a triangle, we have $\ci_{\alpha} = \sum\limits_{1 \leq i < j \leq n} |C(i) - C(j)|^{\alpha} = nk - k^2$. If $n$ is even, then the maximum is attained when $k = n/2$, which yields $n^2/4$. \textit{Thus we have found a graph for which the degree index is zero, but the clustering index is $n^2/4$, the maximal possible value}. 
\begin{figure}[H]
\centering
    \begin{tikzpicture}[scale=0.5]
  \newcommand{\drawPolygon}[3]{
    \begin{scope}[xshift=#1cm, yshift=#2cm]
      \foreach \i in {1,...,#3} {
        \coordinate (v\i) at (360/#3*\i:1.5);
        \fill[red] (v\i) circle(2pt);
      }
      \draw[blue, thick] \foreach \i in {1,...,#3} { (v\i) -- ++(360/#3:1.5) } -- cycle;
    \end{scope}
  }

  \begin{scope}[yshift=7cm]
    \draw[blue, thick] (0,0) -- (1,2) -- (2,0) -- cycle;
    \foreach \i in {0,2} {
      \fill[red] (\i,0) circle(2pt);
    }
    \fill[red] (1,2) circle(2pt);

    \begin{scope}[xshift=5cm, yshift=1cm]
      \foreach \i in {1,2,3,4,5,6,7,8} {
        \coordinate (v\i) at (22.5+360/8*\i:1.5);
        \fill[red] (v\i) circle(2pt);
      }
      \draw[blue, thick] (v1) -- (v2) -- (v3) -- (v4) -- (v5) -- (v6) -- (v7) -- (v8) -- cycle;
    \end{scope}

    \begin{scope}[xshift=8cm]
      \draw[blue, thick] (0,0) -- (1,2) -- (2,0) -- cycle;
      \foreach \i in {0,2} {
        \fill[red] (\i,0) circle(2pt);
      }
      \fill[red] (1,2) circle(2pt);
    \end{scope}
  \end{scope}

  \begin{scope}[yshift=3.5cm]
    \draw[blue, thick] (0,0) -- (1,2) -- (2,0) -- cycle;
    \foreach \i in {0,2} {
      \fill[red] (\i,0) circle(2pt);
    }
    \fill[red] (1,2) circle(2pt);

    \begin{scope}[xshift=5cm, yshift=1cm]
      \foreach \i in {1,2,3,4,5,6} {
        \coordinate (v\i) at (360/6*\i:1.5);
        \fill[red] (v\i) circle(2pt);
      }
      \draw[blue, thick] (v1) -- (v2) -- (v3) -- (v4) -- (v5) -- (v6) -- cycle;
    \end{scope}

    \begin{scope}[xshift=8cm]
      \draw[blue, thick] (0,0) -- (1,2) -- (2,0) -- cycle;
      \foreach \i in {0,2} {
        \fill[red] (\i,0) circle(2pt);
      }
      \fill[red] (1,2) circle(2pt);
    \end{scope}
  \end{scope}

  \begin{scope}[yshift=0cm]
    \draw[blue, thick] (0,0) -- (1,2) -- (2,0) -- cycle;
    \foreach \i in {0,2} {
      \fill[red] (\i,0) circle(2pt);
    }
    \fill[red] (1,2) circle(2pt);

    \begin{scope}[xshift=5cm, yshift=1cm]
      \foreach \i in {1,2,3,4} {
        \coordinate (v\i) at (45+360/4*\i:1.5);
        \fill[red] (v\i) circle(2pt);
      }
      \draw[blue, thick] (v1) -- (v2) -- (v3) -- (v4) -- cycle;
    \end{scope}

    \begin{scope}[xshift=8cm]
      \draw[blue, thick] (0,0) -- (1,2) -- (2,0) -- cycle;
      \foreach \i in {0,2} {
        \fill[red] (\i,0) circle(2pt);
      }
      \fill[red] (1,2) circle(2pt);
    \end{scope}
  \end{scope}
\end{tikzpicture}
\caption{A graph with $d_i = 2$ for any $i$ (hence $\di_{\alpha} = 0$) but $\ci_{\alpha}(\mathcal{G}) = n^2/4$.}
\end{figure}        
        
        Next, we also discuss an example where the clustering index is zero while the degree index $\sum\limits_{1 \leq i < j \leq n} |d_i - d_j|$ is $\frac{n^2(n/2-3)}{4}$. (Lemma \ref{prop:elt} can be used to show that the degree index can be at most $n^3/4$. So our example here is off by a factor of $1/2$, but still demonstrates that the degree index can be very large while the clustering index is zero.) 
        
        Let $n$ be an integer divisible by 6. Consider a complete graph with $n/2$ vertices, and add $n/6$ disjoint triangles to the graph. Then $C(i) = 1$ for any $i$, and so the clustering index is zero. But on the other hand, $|d_i - d_j| = n/2-3$ if one of the $i,j$ is in our original complete graph while the other belongs to a triangle. Since there are $n^2/4$ of these non-zero terms, we conclude  $\sum\limits_{1 \leq i \leq j \leq n} |d_i - d_j| = \frac{n^2(n/2-3)}{4}$.

\begin{figure}[H]
    \centering
    \begin{tikzpicture}
  \draw[blue, thick] (0,0) -- (1,2) -- (2,0) -- cycle;
  \foreach \i in {0,2} {
    \fill[red] (\i,0) circle(2pt);
  }
  \fill[red] (1,2) circle(2pt);

  \begin{scope}[xshift=5cm, yshift=1cm]
    \foreach \i in {1,2,...,6} {
      \coordinate (v\i) at (60*\i:1.5);
      \fill[red] (v\i) circle(2pt);
    }
    \foreach \i in {1,2,...,6} {
      \foreach \j in {\i,...,6} {
        \draw[blue, thick] (v\i) -- (v\j);
      }
    }
  \end{scope}

  \begin{scope}[xshift=8cm]
    \draw[blue, thick] (0,0) -- (1,2) -- (2,0) -- cycle;
    \foreach \i in {0,2} {
      \fill[red] (\i,0) circle(2pt);
    }
    \fill[red] (1,2) circle(2pt);
  \end{scope}
\end{tikzpicture}
\caption{A graph with $C(i) = 1$ for any $i$ (hence $\ci_{\alpha} (\mathcal{G})= 0$) but $\di_{1}(\mathcal{G}) = \frac{n^2(n/2-3)}{4}$.}
\end{figure}

        \end{example} 

        \section{Clustering index of  Erd\H{o}s-R\'enyi  graphs }\label{sec:clusind}

\subsection{Sublinearity of $\mathbb{E}[\ci_1(\mathcal{G})]$ and $\mathbb{E}[\ci_2(\mathcal{G})]$}

In this subsection we will focus on the case where $\mathcal{G}$ is an Erd\H{o}s-R\'enyi graph, and show that $\mathbb{E}[\ci_1(\mathcal{G})]$ and $\mathbb{E}[\ci_2(\mathcal{G})]$ are both sublinear. The result for the latter will later be improved, and in particular $\mathbb{E}[\ci_2(\mathcal{G})]$ will be shown to be bounded by a constant independent of $n$. We begin with the following proposition concerning the first two moments of the local clustering coefficient.

  \begin{proposition}\label{propn:CiCi2}
            Let $\mathcal{G} = (V, E)$ be an Erd\H{o}s-R\'enyi graph with parameters $n \in \mathbb{N}$ and $p \in (0,1)$. Let $i \in V$ be any node. 
            
           (i) 
            We have, 
        $$\mathbb{E}[C(i)] =  p (1 -  (1-p)^{n-1} - (n-1)p(1-p)^{n-2}).$$ In particular, $\mathbb{E}[C(i)] \geq p - \frac{D}{n^2}$ for each $n$, for some constant $D$  depending on $p$, but   not on $n$. 

(ii) We have, 
            $$
         \left| \mathbb{E}[C(i)^2] -   p^2 \right| \leq \frac{D_1}{n^2},$$ where $D_1$ is a positive constant depending on $p$, but not on $n$. 
        
       (iii) $\var(C(i)) = \mathcal{O}\left(\frac{1}{n^2}\right).$
        \end{proposition}

\textbf{Proof.} (i) First, letting $\phi = \mathbb{P}(d_i \leq 1) = \mathbb{P}(d_i = 0) + \mathbb{P}(d_i = 1) = (1-p)^{n-1} + (n-1)p(1-p)^{n-2}$, we have 
\begin{equation}\label{eqn:eltc}
\mathbb{E}[C(i)] = \phi   \mathbb{E}[C(i) \mid d_i \leq 1] + (1 - \phi)   \mathbb{E}[C(i) \mid d_i \geq 2].     
\end{equation}
    Now, by the definition of $C(i)$, $ \mathbb{E}[C(i) \mid d_i \leq 1] = 0$. We claim  $ \mathbb{E}[C(i) \mid d_i \geq 2] = p$.  We have 
$$\mathbb{E}[C(i) \mid d_i \geq 2] = \mathbb{E}[ \mathbb{E}[ C(i) \mid d_i, d_i \geq 2 ] \mid d_i \geq 2] = \mathbb{E} \left[ \frac{1}{\binom{d_i}{2}} \mathbb{E}  \left[ \sum_{k, \ell \in N(i)} \mathbf{1}(k \leftrightarrow \ell) \mid d_i, d_i \geq 2   \right] \mid d_i \geq 2 \right],$$ where we recall that $N(i)$ is the set of neighbors of node $i$. Now given $d_i$ with $d_i \geq 2$, $\sum_{k, \ell \in N(i)} \mathbf{1}(k \leftrightarrow \ell)$ is binomially distributed with parameters $\binom{d_i}{2}$ and $p$, whose expectation is $\binom{d_i}{2} p$. Therefore, continuing the last observations, we conclude, $$\mathbb{E}[C(i) \mid d_i \geq 2]  = \mathbb{E} \left[ \frac{1}{\binom{d_i}{2}} \binom{d_i}{2} p \right] = p.$$
Substituting this into \eqref{eqn:eltc}, we obtain 
$$\mathbb{E}[C(i)] =  p (1 -  (1-p)^{n-1} - (n-1)p(1-p)^{n-2}) = p - p(1-p)^{n-1}  - (n-1)p^2 (1-p)^{n-2}.$$ Now, for $p \in (0,1)$, the sequence  $p(1-p)^{n-1} + (n-1)p^2 (1-p)^{n-2}$ decays exponentially fast as $n \rightarrow \infty$  and thus, it is  upper bounded by $D / n^2$ for some $D > 0$.
The assertion that $\mathbb{E}[C(i)] \geq p - \frac{D}{n^2}$ is now clear.

(ii)
Note that below, $D$ denotes a constant independent of $n$ that does not necessarily have the same value in its two appearances. 

        Let $A_i$ be the event that $|d_i - (n - 1)p| <  n^{2/3}$. Then by McDiarmid (or Azuma-Hoeffding) inequality \cite{mcd:98}, $\mathbb{P}(\overline{A_i}) \leq M_1 e^{-M_2 n^{1/3}}$ for some positive constants $M_1, M_2$ independent of $n$. Using $A_i$, we now write, 
$$\mathbb{E}[C(i)^2] = \mathbb{E}[C(i)^2 \mathbf{1}(A_i)] + \mathbb{E}[C(i)^2 \mathbf{1}(\overline{A_i})].$$ 
Let us analyze the two terms on the right-hand side separately. First, focusing on the second term, noting the trivial bound $|C(i)| \leq 1$, observe that $$\mathbb{E}[C(i)^2 \mathbf{1}(\overline{A_i})] \leq \mathbb{P}(\overline{A_i}) \leq  M_1 e^{-M_2 n^{1/3}} \leq \frac{D}{n^2}.$$
Hence for the second term, we have, 
\begin{eqnarray}\label{est:2nd}
    0 \leq \mathbb{E}[C(i)^2 \mathbf{1}(\overline{A_i})]  \leq \frac{D}{n^2}.
\end{eqnarray}
Moving on to the first term, we have
\begin{eqnarray}\label{eqn:useful}
   \nonumber \mathbb{E}[C(i)^2 \mathbf{1}(A_i)] &=& \mathbb{E}[\mathbb{E} [C(i)^2 \mathbf{1}(A_i) \mid d_i ]] \\
\nonumber    &=& 
    \mathbb{E}\left[\mathbb{E} \left[\frac{1}{\binom{d_i}{2}^2} \left( \sum_{\{j,k\} \subset N(i)} \mathbf{1}(j \leftrightarrow k)\right)^2 \mathbf{1}(A_i) \mid d_i \right]\right] \\
\nonumber    &=& 
    \mathbb{E}\left[\frac{\mathbf{1}(A_i)}{\binom{d_i}{2}^2} \mathbb{E} \left[ \sum_{\{j,k\} \subset N(i)} \mathbf{1}(j \leftrightarrow k) + \sum_{\{j,k\}, \{j',k'\} \subset N(i),  \{j,k\} \neq \{j',k'\}} \mathbf{1}(j \leftrightarrow k, j' \leftrightarrow k') \mid d_i \right]\right] \\
 \nonumber   &=& 
    \mathbb{E}\left[\frac{\mathbf{1}(A_i)}{\binom{d_i}{2}^2}  \left(  \binom{d_i}{2} p + \left( \binom{d_i}{2}^2 - \binom{d_i}{2} \right) p^2 \right) \right] \\
    &=& p^2 \mathbb{E}[\mathbf{1}(A_i)] + (p - p^2) \mathbb{E} \left[ \frac{\mathbf{1}(A_i)}{\binom{d_i}{2}} \right] = p^2 \mathbb{P}(A_i) + (p - p^2) \mathbb{E} \left[ \frac{\mathbf{1}(A_i)}{\binom{d_i}{2}} \right]. 
\end{eqnarray}
Now, when $A_i$ is true, $(n - 1)p - n^{2 / 3} \leq d_i \leq (n - 1)p + n^{2 / 3}$. Keeping this in mind, and using the last relation we obtained, we get, 
$$ \mathbb{E}[C(i)^2 \mathbf{1}(A_i)] \leq p^2 \cdot 1 + \frac{p-p^2}{\binom{(n - 1)p - n^{2 / 3}}{2}} \leq p^2 + \frac{D}{n^2}.$$
Also, again by \eqref{eqn:useful}, we have $$ \mathbb{E}[C(i)^2 \mathbf{1}(A_i)]  \geq  p^2 \mathbb{P}(A_i) \geq p^2 (1 - M_1 e^{-M_2 n^{1/3}})  \geq p^2 - \frac{D}{n^2}.$$ 
Hence, we conclude, 
\begin{equation}\label{est:1st}
    |\mathbb{E}[C(i)^2 \mathbf{1}(A_i)] - p^2 | \leq \frac{D}{n^2}. 
\end{equation}

Combining \eqref{est:2nd} and \eqref{est:1st}, we arrive at, 
$$|\mathbb{E}[C(i)^2 - p^2 ] | \leq \frac{D_1}{n^2}, $$ which holds for every $n$, where $D_1$ is depending on $p$, but not on $n$, as asserted.

(iii) Follows immediately from (i) and (ii). \hfill $\square$

Proposition \ref{propn:CiCi2} can now be used to show that  $\mathbb{E}[\ci_1(\mathcal{G})]$ is  sublinear.
         
        \begin{theorem}\label{thm:ci1}
            Let $\mathcal{G}$ be an Erd\H{o}s-R\'enyi graph with parameters $n \in \mathbb{N}$ and $p \in (0,1)$. Then we have $$\mathbb{E}[\ci_1(\mathcal{G})] \leq K n, \quad n \geq 1,$$  where $K$ is a constant depending on $p$, but not on $n$.

        \end{theorem}

        \textbf{Proof.} Let $i, j \in \mathbb{N}$. Note that $C(i)$ and $C(j)$ have the same distribution in Erd\H{o}s-R\'enyi graphs due to the underlying symmetry. Using this along with the Cauchy-Schwarz inequality, and the previous proposition, we observe,
        \begin{eqnarray*}
        \mathbb{E}[|C(i) - C(j)|] &=& \mathbb{E}[|C(i)  - \mathbb{E}[C(i)] + \mathbb{E}[C(j)]- C(j)|] \\
        &\leq& 2 \mathbb{E} |C(i) - \mathbb{E}|C(i)||
        \leq 2 (\mathbb{E} |C(i) - \mathbb{E}|C(i)||^2)^{1/2} 
        = 2 \sqrt{\var(C(i))} \leq \frac{D}{n}, 
        \end{eqnarray*}
        for some constant $D$. 
        Thus, 
        $$
        \mathbb{E}[\ci_1(\mathcal{G})]  \leq \binom{n}{2} \frac{D}{n} \leq K n. $$
         \hfill $\square$

Noting that the trivial inequality $|C(i) - C(j)|^2 \leq |C(i) - C(j)|$ always holds, one obtains the following corollary. 

\begin{corollary}\label{Cor:ci2}
    In the setting of the previous theorem, $$\mathbb{E}[\ci_2(\mathcal{G})] \leq K n, \quad n \geq 1,$$  where $K$ is a constant depending on $p$, but independent of $n$.  
\end{corollary}

In the next subsection, we will show that in the case of Erd\H{o}s-R\'enyi graphs for fixed $p$, $\mathbb{E}[\ci_2(\mathcal{G})]$ is indeed bounded by a constant that does not depend on $n$. Note that our simulation results further suggest that for Erd\H{o}s-R\'enyi graphs, the growth of $\ci_1(\mathcal{G})$ behaves linearly in the number of nodes, and $\ci_2(\mathcal{G})$ converges to a positive constant as the number of nodes increases. 

\subsection{Further analysis of $\mathbb{E}[\ci_2(\mathcal{G})]$}

        Now that we know $\mathbb{E}[\ci_2(\mathcal{G})]$ is sublinear via Corollary \ref{Cor:ci2} when $\mathcal{G}$ is an Erd\H{o}s-R\'enyi graph, we will focus on this expectation in more detail and provide a more precise analysis for this case. This will require understanding of expectations of the form $\mathbb{E}[C(i) C(j)]$ for nodes $i \neq j$. 

        \begin{proposition}\label{propn:poscor}
            Let $\mathcal{G}$ be an Erd\H{o}s-R\'enyi graph with parameters $n \in \mathbb{N}$ and $p \in (0,1)$, and $i,j $ be two distinct nodes. Then, $$\mathbb{E}[C(i) C(j)] \geq p^2 - \frac{D_2}{n^2}, $$  for each $n$, where $D_2$ is a positive constant depending on $p$, but not on $n$.
        \end{proposition}

        We defer the proof of Proposition \ref{propn:poscor} to the end of this subsection and first discuss the main result. 

           \begin{theorem}
            Let $\mathcal{G}$ be an Erd\H{o}s-R\'enyi graph with parameters $n \in \mathbb{N}$ and $p \in (0,1)$. Then we have $$\mathbb{E}[\ci_2(\mathcal{G})] \leq K $$ for some constant $K$ independent of $n$.
        \end{theorem}
        
        \textbf{Proof.} By Proposition \ref{propn:CiCi2} and Proposition \ref{propn:poscor}, we have:      
        \begin{eqnarray*}   
            \mathbb{E}[\ci_2(\mathcal{G})] &=& \binom{n}{2} \mathbb{E}[(C(i)-C(j))^2] = \binom{n}{2} \left( \mathbb{E}[C(i)^2] + \mathbb{E}[C(j)^2] - 2\mathbb{E}[C(i)C(j)] \right) \\
            && = n(n-1) \left( \mathbb{E}[C(i)^2] - \mathbb{E}[C(i)C(j)] \right) \leq n^2 \left(p^2 + \frac{D_1}{n^2} - p^2 + \frac{D_2}{n^2} \right) = D_1 + D_2 =:K.
        \end{eqnarray*}
        \hfill $\square$

  Let us now prove Proposition \ref{propn:poscor}.

        \textbf{Proof of Proposition \ref{propn:poscor}.} 
       Recall from the proof of Proposition \ref{propn:CiCi2} that  $A_i$ is the event that $|d_i - (n - 1)p| <  n^{2/3}$. Define $A_j$ similarly. 
 We have $$\mathbb{E}[C(i)C(j) ] = \mathbb{E}[C(i)C(j) \mathbf{1}(A_i) \mathbf{1}(A_j) ] + S,$$ where $$S=  \mathbb{E}[C(i)C(j) \mathbf{1}(A_i) \mathbf{1}(\overline{A_j}) ] + \mathbb{E}[C(i)C(j) \mathbf{1}(\overline{A_i}) \mathbf{1}(A_j) ]  + \mathbb{E}[C(i)C(j) \mathbf{1}(\overline{A_i}) \mathbf{1}(\overline{A_j}) ].$$ Noting that $A_i$ and $A_j$ are equally likely, and  using trivial upper bound for the clustering coefficients, we see that $$S \leq 3 \mathbb{P}(\overline{A_i}).$$ Recalling now $\mathbb{P}(\overline{A_i}) \leq M_1 e^{-M_2 n^{1/3}}$ for some constants $M_1, M_2 > 0$  from the Proof of Proposition 3.1, we conclude  that $$S\leq \frac{D}{n^2}, \quad n \geq 1,$$ where $D$ is a positive constant independent of $n$. (Again, $D$'s in distinct appearances may denote different constants.) Next, we focus on the estimation of $\mathbb{E}[C(i)C(j) \mathbf{1}(A_i) \mathbf{1}(A_j) ] $. 
 
Let $M$ be the event that the edge between $i, j$ exists. Also, let $T_i$ and $T_j$ denote the number of triangles that contain $i$ and $j$ as a node, respectively. Then 
   \begin{eqnarray}\label{eqn:mainest1}
  \nonumber \mathbb{E}[C(i)C(j) \mathbf{1}(A_i) \mathbf{1}(A_j) ]  &=& p \mathbb{E}[C(i)C(j) \mathbf{1}(A_i) \mathbf{1}(A_j) \mid M] + (1 - p) \mathbb{E}[C(i)C(j) \mathbf{1}(A_i) \mathbf{1}(A_j) \mid \overline{M}]   \\
  &=& p \mathbb{E} \left[ \frac{\mathbf{1}(A_i) \mathbf{1}(A_j)}{\binom{d_i}{2} \binom{d_j}{2}} \mathbb{E}[T_i T_j \mid d_i, d_j, M] \mid M\right] \\
 \nonumber && + (1- p) \mathbb{E} \left[ \frac{\mathbf{1}(A_i) \mathbf{1}(A_j)}{\binom{d_i}{2} \binom{d_j}{2}} \mathbb{E}[T_i T_j \mid d_i, d_j, \overline{M}] \mid \overline{M}\right]. 
\end{eqnarray}

Now we will examine the conditional expectation of $T_i T_j$ given $d_i, d_j, M$ and again of $T_i T_j$ given $d_i, d_j, \overline{M}$ in more detail. For this purpose, we decompose $T_i$ and $T_j$ as $$ T_i = T_{i1} + T_{i2} + T_{i3} \quad \text{and} \quad T_j = T_{j1} + T_{j2} + T_{j3},$$ where, 
\begin{itemize}
    \item $T_{i1} = T_{j1}$ denotes the number of triangles that have both $i$ and $j$ as vertices,
    \item $T_{i2} = T_{j2}$ denotes the number of edges $e_{k \ell}$, with $k$ and $\ell$ distinct from $i$ and $j$, and the triangles $\{i,k,\ell\}$ and $\{j,k,\ell\}$ both appear in the graph, 
    \item $T_{i3}$ and $T_{j3}$ denote the remaining triangles containing $i$ and $j$ as a node, respectively. 
\end{itemize}

A few observations about this decomposition are in order. 
First note that conditional on $\overline{M}$, both $T_{i1}$ and $T_{j1} = 0.$ Second, apart from the products $T_{i1}T_{j1} = T_{i1}^2$ and $T_{i2}T_{j2} = T_{i2}^2$, the products appearing in $T_i T_j$ are products of independent random variables given $M, d_i, d_j$ or $\overline{M}, d_i, d_j$.
Hence, to estimate these terms, we can investigate the remainder we obtain from $\mathbb{E}[T_i]\mathbb{E}[T_j]$ in the following way.  
We have,  
\begin{align}\label{eqn:mainest2}
  \nonumber  \mathbb{E}[T_i T_j \mid M, d_i, d_j] &= \sum_{k=1}^3 \sum_{\ell =1}^3 \mathbb{E}[T_{ik}T_{j\ell} \mid M,d_i,d_j] \\
 \nonumber   &= \left(\sum_{(k, \ell) \notin \{(1,1),(2,2)\}} \mathbb{E}[T_{ik}  \mid M,d_i,d_j] \mathbb{E}[T_{j\ell}  \mid M,d_i,d_j]  \right)\\
    \nonumber & \quad + \mathbb{E}[T_{i1}^2  \mid M,d_i,d_j]  + \mathbb{E}[T_{i2}^2  \mid M,d_i,d_j]\\
   \nonumber &= \sum_{k=1}^3 \sum_{\ell =1}^3  \mathbb{E}[T_{ik}  \mid M,d_i,d_j] \mathbb{E}[T_{j\ell}  \mid M,d_i,d_j]   \\
   \nonumber & \qquad -  \mathbb{E}[T_{i1}  \mid M,d_i,d_j] \mathbb{E}[T_{j1}  \mid M,d_i,d_j] -  \mathbb{E}[T_{i2}  \mid M,d_i,d_j] \mathbb{E}[T_{j2}  \mid M,d_i,d_j] \\
  \nonumber  & \qquad  +  \mathbb{E}[T_{i1}^2  \mid M,d_i,d_j] + \mathbb{E}[T_{i2}^2  \mid M,d_i,d_j] \\
  \nonumber  &=  \mathbb{E}[T_{i}  \mid M,d_i,d_j] \mathbb{E}[T_{j}  \mid M,d_i,d_j]  \\
  \nonumber  & \qquad +\left( \mathbb{E}[T_{i1}^2  \mid M,d_i,d_j] - \left(\mathbb{E}[T_{i1}  \mid M,d_i,d_j] \right)^2 \right) \\
  \nonumber  & \qquad + \left( \mathbb{E}[T_{i2}^2  \mid M,d_i,d_j] - \left(\mathbb{E}[T_{i2}  \mid M,d_i,d_j] \right)^2 \right)   \\
   &= \mathbb{E}[T_i \mid M, d_i, d_j] \mathbb{E}[T_j \mid M, d_i, d_j]\\
   \nonumber & \qquad  + \mathbb{E}[ T_{i1}^2 - (\mathbb{E}[T_{i1}|M,d_i,d_j])^2|M,d_i, d_j]  + \mathbb{E}[ T_{i2}^2 - (\mathbb{E}[T_{i2}|M,d_i,d_j])^2|M,d_i, d_j].
\end{align}
Also, recalling that conditional on $\overline{M}$, both $T_{i1}$ and $T_{j1} = 0$, and following similar  steps we obtain, 
\begin{equation}\label{eqn:mainest3}
    \mathbb{E}[T_i T_j \mid \overline{M}, d_i, d_j] = \mathbb{E}[T_i \mid \overline{M}, d_i, d_j] \mathbb{E}[T_j \mid \overline{M}, d_i, d_j] + \mathbb{E}[ T_{i2}^2 - (\mathbb{E}[T_{i2}|\overline{M},d_i,d_j])^2|\overline{M}, d_i, d_j].
\end{equation}
Now, note that 
\begin{equation}\label{eqn:mainest4}
\mathbb{E}[T_i \mid  M, d_i, d_j] = \mathbb{E}[T_i \mid \overline{M}, d_i, d_j] = p \binom{d_i}{2}.    
\end{equation}
  To see this, given the degree $d_i$ of $i$th node along with $d_j$ and $M$, observe that there are $\binom{d_i}{2}$ possible triangles containing node $i$, and one just needs to count the number of edges between the pairs of neighbors of $i$ to find $T_i$. Hence, there are $\binom{d_i}{2}$ candidates for triangles containing node $i$, each of which is present with probability $p$, independent of others. Thus, one has $$\mathbb{E}[T_i \mid  M, d_i, d_j]  = p \binom{d_i}{2}.$$ Noting that given the degree of $d_i$, the presence or absence of the edge between node $i$ and $j$ does not change the distribution of $T_i$, we similarly have $$\mathbb{E}[T_i \mid \overline{M}, d_i, d_j] = p \binom{d_i}{2}.$$

Combining our observations in \eqref{eqn:mainest1}, \eqref{eqn:mainest2}, \eqref{eqn:mainest3} and \eqref{eqn:mainest4}, we obtain,
\begin{align*}
    \mathbb{E}[C(i)C(j) \mathbf{1}(A_i) \mathbf{1}(A_j) ] &= p\mathbb{E} \left[ \frac{\mathbf{1}(A_i) \mathbf{1}(A_j)}{\binom{d_i}{2} \binom{d_j}{2}}  p^2\binom{d_i}{2} \binom{d_j}{2}  \mid M \right] \\
    &+ (1-p)\mathbb{E} \left[ \frac{\mathbf{1}(A_i) \mathbf{1}(A_j)}{\binom{d_i}{2} \binom{d_j}{2}}  p^2\binom{d_i}{2} \binom{d_j}{2}  \mid \overline{M} \right] \\
    &+ p \mathbb{E} \left[ \frac{\mathbf{1}(A_i) \mathbf{1}(A_j)}{\binom{d_i}{2} \binom{d_j}{2}} \mathbb{E}[ T_{i1}^2 - (\mathbb{E}[T_{i1}|M,d_i,d_j])^2|M,d_i, d_j] \mid M\right] \\
    &+ p \mathbb{E} \left[ \frac{\mathbf{1}(A_i) \mathbf{1}(A_j)}{\binom{d_i}{2} \binom{d_j}{2}} \mathbb{E}[ T_{i2}^2 - (\mathbb{E}[T_{i2}|M,d_i,d_j])^2|M,d_i, d_j] \mid M\right] \\
    &+ (1- p) \mathbb{E} \left[ \frac{\mathbf{1}(A_i) \mathbf{1}(A_j)}{\binom{d_i}{2} \binom{d_j}{2}} \mathbb{E}[ T_{i2}^2 - (\mathbb{E}[T_{i2}|\overline{M},d_i,d_j])^2|\overline{M},d_i, d_j] \mid \overline{M}\right].
\end{align*}
The first two terms on the right-hand side simplify to $$p^2 \left(p \mathbb{E} \left[ \mathbf{1}(A_i) \mathbf{1}(A_j) \mid M \right] + (1-p) \mathbb{E} \left[ \mathbf{1}(A_i) \mathbf{1}(A_j) \mid \overline{M} \right] \right),$$ which is equal to $p^2\mathbb{E} \left[ \mathbf{1}(A_i) \mathbf{1}(A_j) \right]$. Also, all the other terms contain conditional variances of $T_{i1}$ and $T_{i2}$, which are bound to be non-negative random variables and hence, these terms are non-negative. 
Thus, we conclude that 
$$ \mathbb{E}[C(i)C(j) \mathbf{1}(A_i) \mathbf{1}(A_j) ]  \geq p^2\mathbb{E} \left[ \mathbf{1}(A_i) \mathbf{1}(A_j) \right].$$
Also, again recalling from the proof of Proposition \ref{propn:CiCi2}  that $\mathbb{P}(\overline{A_i}) \leq M_1 e^{-M_2 n^{1/3}}$  for some constants $M_1, M_2$ independent of $n$, and noting a similar bound holds for $A_j$, we have $\mathbb{P}(A_i \cap A_j) \geq 1 - \frac{D}{n^2}$ for some constant $D$ independent of $n$. Using this, we reach at, 
$$ \mathbb{E}[C(i)C(j) \mathbf{1}(A_i) \mathbf{1}(A_j) ] \geq p^2 - \frac{D}{n^2},$$ where $D$ is again independent of $n$. 

 Lastly, recalling $$\mathbb{E}[C(i)C(j) ] = \mathbb{E}[C(i)C(j) \mathbf{1}(A_i) \mathbf{1}(A_j) ] + S,  $$ where $0 \leq S \leq \frac{D}{n^2}$, we get, 
 $$\mathbb{E}[C(i)C(j) ] \geq p^2 - \frac{D_2}{n^2},$$ as asserted. \hfill $\square$

         \section{Degree index in Erd\H{o}s-R\'enyi graphs}\label{sec:degindex}

        As already noted earlier, the degree index (or total degree irregularity) has been previously introduced and studied, but we are not aware of a relevant analysis for random graphs. This section is devoted to such an analysis, and this turns out to be more tractable compared to the clustering index case due to the lack of $\binom{d_i}{2}$ terms in the denominator.  
         We begin the discussion by studying $\mathbb{E}[\di_2(\mathcal{G})]$. 

        \begin{theorem}\label{thm:di2}
        For an Erd\H{o}s-R\'enyi graph $\mathcal{G}$ with parameters $n \in \mathbb{N}$ and $p \in (0,1)$, $$\mathbb{E}[\di_2(\mathcal{G})] = 6\binom{n}{3}p(1-p).$$ In particular, 
         $$\frac{\mathbb{E}[\di_2(\mathcal{G})]}{6\binom{n}{3}} = p(1-p),$$ which is independent of $n$. 
        \end{theorem}

        \textbf{Proof. }
        Let $\mathcal{G}$ be an Erd\H{o}s-R\'enyi graph with parameters $n \in \mathbb{N}$ and $p \in (0,1)$.  We are interested in calculating $$\mathbb{E}[(d_i - d_j)^2] = 2\mathbb{E}[d_i^2] - 2\mathbb{E}[d_i d_j].$$
         Let $k = n-2$  for convenience. Denoting the binomial distribution by $\bin$, since $d_i \sim \bin(n-1,p) = \bin(k+1,p)$, we have 
        $$\mathbb{E}[d_i^2] = (\mathbb{E}[d_i])^2 + \var(d_i) = (n-1)^2p^2 + (n-1)p(1-p) = (k+1)^2p^2 + (k+1)p(1-p).$$
        Next, let us examine $\mathbb{E}[d_i d_j]$. Let $M$ denote the event that $e_{ij}$ appears in the graph. Conditioning on $M$ and $\overline{M}$ yields 
        \begin{align*}
            \mathbb{E}[d_i d_j] = p\mathbb{E}[d_id_j \mid M] + (1-p) \mathbb{E}[d_id_j \mid \overline{M}].
        \end{align*}
        
        If we assume that $e_{ij}$ does not appear in the graph, then the number of possible neighbors of $i$ or $j$ is now $n-2$ instead of $n-1$. And since we are in the Erd\H{o}s-R\'enyi model, the presence of any other edge is independent of the presence of $e_{ij}$. Moreover, apart from whether $e_{ij}$ appears or not, $d_i$ and $d_j$ are independent. By the above observations, $d_i \mid M$ and $d_j  \mid  M$ are i.i.d. with the common distribution $1 + \bin(n-2,p) = 1 + \bin(k,p)$. Similarly, $d_i \mid \overline{M}$ and $d_j  \mid \overline{M}$ are i.i.d. this time with the common distribution $\bin(n-2,p) = \bin(k,p)$. Therefore, we have
        \begin{align*}
            \mathbb{E}[d_i d_j] &= p\mathbb{E}[d_i d_j  \mid M] + (1-p)\mathbb{E}[d_i d_j  \mid \overline{M}] = p (kp + 1)^2 + (1-p) (kp)^2 \\ &= k^2p^3 + 2kp^2 + p + k^2p^2 - k^2p^3 = k^2p^2 + 2kp^2 + p,
        \end{align*}
        which implies
        \begin{eqnarray*}
            \mathbb{E}[d_i^2] - \mathbb{E}[d_i d_j] &=& (k+1)^2 p^2 + (k+1)p(1-p) - k^2p^2 - 2kp^2 - p \\
            &=& k^2p^2 + 2kp^2 + p^2 + kp - kp^2 + p - p^2 - k^2p^2 - 2kp^2 - p \\
            &=& kp - kp^2 = kp(1-p) = (n-2)p(1-p).
        \end{eqnarray*}
        Thus, $$\mathbb{E}[(d_i-d_j)^2] = 2(n-2)p(1-p).$$ Recalling that $\di_{2} (\mathcal{G} )= \sum_{1 \leq i < j \leq n} |d_i - d_j|^{2}$,  result follows.  \hfill $\square$

        Next, we focus on $\di_1(\mathcal{G})$. In this case a direct computation readily gives
        $$\mathbb{E}[\di_1(\mathcal{G})] = \mathbb{E} \left[\sum_{1 \leq i < j \leq n} |d_i - d_j| \right] = \sum\limits_{\ell =0}^{n-2} \sum\limits_{k=0}^{n-2} |k-\ell| \binom{n-1}{k} \binom{n-1}{\ell} p^k (1-p)^{n-1-k} p^{\ell} (1-p)^{n-1-\ell}.$$
        For clearer expressions, we will analyze the symmetric $p=1/2$ case and the general $p$ case separately. 
        
        Keeping the notations above, let us now observe that 
        \begin{align*}
            \mathbb{E}[|d_i - d_j|] = p\mathbb{E}[|d_i - d_j| \mid M] + (1-p) \mathbb{E}[|d_i - d_j| \mid \overline{M}].
        \end{align*} 
        Note that the distribution of $d_i - d_j$ is the same regardless of whether the edge between $i$ and $j$, $e_{ij}$, is present in the graph or not. And since the only dependence between $d_i$ and $d_j$ was a result of $e_{ij}$, we have, $$\mathbb{E}[|d_i - d_j|] = \mathbb{E}[|B_1 - B_2|],$$ where $B_1, B_2$ are independent binomial random variables with parameters $n - 2$ and $p$.  
        For general $p$ case, an easy upper bound on $\mathbb{E}[|d_i - d_j|]$ can be obtained as follows:   
        \begin{eqnarray*}
            \mathbb{E}[|d_i - d_j|] &=& \mathbb{E}|B_1 - B_2| \leq \mathbb{E}|B_1 - (n-2)p| + \mathbb{E}|B_2 - (n-2)p| = 2 \mathbb{E}|B_1 - (n-2)p|  \\ && \leq 2 (\mathbb{E}|B_1 - (n-2)p|^2 )^{1/2} = 2 \sqrt{\var(B_1)} = 2 \sqrt{(n-2) p (1-p)}.
        \end{eqnarray*}
        Improving this, one may obtain an exact asymptotic result for $\mathbb{E}[\di_1(\mathcal{G})]$. This will follow from the following elementary result on the mean absolute difference of binomials. 
        
        \begin{lemma}\label{propn:bin}
            Let $X_n, Y_n$ be independent binomial random variables with parameters $n \in \mathbb{N}$ and $p \in (0,1)$.
        
            (i) If $p = 1/2$, then $$\mathbb{E}|X_n - Y_n|= \frac{n \binom{2n}{n}}{2^{2n}}.$$ Moreover, $\mathbb{E}|X_n - Y_n| \sim \sqrt{\frac{n}{\pi}}$, as $n \rightarrow \infty$. 
        
            (ii) For general $p \in (0,1)$, we have $$\left| \mathbb{E}|X_n - Y_n| -  \frac{2}{\sqrt{\pi}} \sqrt{n p (1-p)}  \right| \leq D$$ for every $n \geq 1$, where $D$ is a constant depending on $p$, but not on $n$. In particular, $\mathbb{E}|X_n - Y_n|  \sim   \frac{2}{\sqrt{\pi}}\sqrt{n p (1-p)}$, as $n \rightarrow \infty$. 
        \end{lemma}

        Although the lemma is elementary and there are similar results in the literature (see, for example, \cite{dz}), we were unable to find the exact form used in our case. Therefore, we include the proof of the lemma in the appendix.

Now, given Lemma \ref{propn:bin}, the following result on $\di_1$  follows immediately. 
        
        \begin{theorem}\label{thm:di1}
         (i) If $\mathcal{G}$  is an Erd\H{o}s-R\'enyi graph with parameters $n \in \mathbb{N}$ and $p = 1/2$, then, 
        $$\mathbb{E}[\di_1(\mathcal{G})] =  \frac{(n-2) \binom{n}{2} \binom{2n-4}{n-2}}{2^{2n-4}}.$$ 
          In particular,  $\mathbb{E}[\di_1(\mathcal{G})] \sim \sqrt{\frac{n-2}{\pi}} \binom{n}{2}$, as $n \rightarrow \infty$.

          (ii)      If $\mathcal{G}$  is an Erd\H{o}s-R\'enyi graph with parameters $n \in \mathbb{N}$ and $p \in (0,1)$, then, 
        $$\mathbb{E}[\di_1(\mathcal{G})] \sim \frac{2}{\sqrt{\pi}} \binom{n}{2} \sqrt{(n-2) p (1 - p)}$$   as $n \rightarrow \infty$.
        \end{theorem}
        

          \begin{remark}\label{rem:norm}
One may consider the normalized versions of the degree indices we study above in order to be able to do comparisons among distinct random graphs models in a clear way and to make meaningful inferences about the graph properties. Motivated by Theorem \ref{thm:di2} and Theorem \ref{thm:di1} above, a natural candidate for the scaled version of the statistic we study is $$\di_{\alpha}^* = \frac{\di_{\alpha} (\mathcal{G}) }{n^{2 + \alpha /2}}.$$
Note that in the case of Erd\H{o}s-R\'enyi graphs, as $n\rightarrow \infty$, $\mathbb{E}[\di_{\alpha}^*]$ converges to   $\sqrt{\frac{p(1-p)}{\pi}}$ and $p(1-p)$  for $\alpha $ equals 1 and 2, respectively.  

Also, a second natural possible normalization for $\di_{\alpha}(\mathcal{G})$ could be a scaling by $\binom{n}{2}$, where $n$ is the number of nodes. However, such a choice does not take $\alpha$ into account, and the scaled limits do not tend to constant values. 
         \end{remark}

         \section{Monte Carlo simulations for other random graph models}\label{sec:otherrandom} 

         The calculations for the degree and clustering indices become more involved when one leaves the framework of   Erd\H{o}s-R\'enyi random graphs and deals with other models. In this section, we focus on three other random graph models and do a Monte Carlo simulation study for these. The models we will be interested in are the Watts-Strogatz model, Barab\'asi-Albert model and lastly, random regular graphs. The degree index in random regular graphs is already zero, and one is only interested in the clustering index for this special case.

\subsection{Random graph models and implementation}\label{sec:otherrgm}

        Here, instead of providing a detailed technical background on the random graph models of interest, we will briefly discuss the practical aspects of the implementation of Monte Carlo simulations. In particular, we will focus on the selection of parameters so that the experiments are performed with matching edge densities for distinct random graph models. Before moving further, let us note that we used the NetworkX module in Python for all the simulations done below. Now, let us review the models under study, keeping in mind that we are willing to have a given edge density $p^* \in (0,1)$  in our simulations. 

        1. \textit{Erd\H{o}s-R\'enyi model}. The function erdos\textunderscore renyi\textunderscore graph which has parameters $n$ and $p$ as discussed in previous sections produces an Erd\H{o}s-R\'enyi graph in NetworkX. One just chooses $p = p^*$ in order to satisfy the edge density criterion. 
        
        2. \textit{Watts-Strogatz model}. The model was introduced by D. J. Watts and S. Strogatz in 1998, and it catches certain common properties of real-life complex networks by having the small world property, high clustering, and short average path lengths \cite{ws:98}. The description of watts\textunderscore strogatz\textunderscore graph function in NetworkX reference page is as follows: \textit{First create a ring over $n$ nodes. Then each node in the ring is joined to its $k$ nearest neighbors (or $k-1$ neighbors if is odd). Then shortcuts are created by replacing some edges as follows: for each edge $(u,v)$ in the underlying “$n$-ring with $k$ nearest neighbors” with probability $p$ (rewiring probability) replace it with a new edge $(u,w)$ with uniformly random choice of existing node $w$.}

        The number of edges in the Watts-Strogatz model is constant and is given by $\frac{nk}{2}$. Thus, in order to have the edge density $p^*$, one needs to have $p^* = \frac{\frac{nk}{2}}{\binom{n}{2}}$. Simplifying, we obtain $p^* = \frac{k}{n - 1}$. Thus, given $p^*$ and $n$, one needs to choose $k$ properly so that this criterion is satisfied.

         3. \textit{Barab\'asi-Albert model}. The model, which was introduced by Albert and  Barab\'asi in 2002, is a random graph model  \cite{ab:02}
         that uses preferential attachment (i.e. higher-degree vertices tend to have even higher degrees) in order to generate scale-free networks. In NetworkX, samples from this model are generated via the barabasi\textunderscore albert\textunderscore graph(n, m, seed=None, initial\textunderscore graph=None) function.  The descriptions of the inputs in NetworkX references page  are as follows. 
         \begin{itemize}
             \item $n$ int Number of nodes
             \item $m$ int Number of edges to attach from a new node to existing nodes
             \item Initial\textunderscore graph Graph or None (default) Initial network for Barab\'asi–Albert algorithm. It should be a connected graph for most use cases. A copy of initial\textunderscore graph is used. If None, starts from a star graph on $(m+1)$ nodes.
         \end{itemize}
        In our case, we use the default as the initial graph and so begin with a star on $m+1$ nodes. Now, given the edge density $p^*$ and the number of nodes $n$, our goal is to select $m$ so that the resulting Barab\'asi-Albert graph will have approximate edge density $p^*$. 

        Now, since we begin with $m+1$ nodes, the steps required to reach $n$ nodes is $N = n - (m+1)$. Then the number of edges inserted throughout these steps is $Nm = (n  - (m +1))m$. Thus, also considering  the initial nodes present, at the end of the evolution of the graph, we have $$(n- (m+1) )m + m$$ many edges, and so the edge density is $$\frac{(n - (m+1))m + m}{\binom{n}{2}} = \frac{nm - m^2}{\frac{n (n-1)}{2}}.$$ Hence we would like to have $\frac{nm - m^2}{\frac{n (n-1)}{2}} \approx p^*$, or rearranging, $$2m^2 -2n m + n(n-1) p^* \approx 0.$$ Now the corresponding $\Delta$ is $4n^2 - (4)  (2) n(n-1)p^* = 4n^2 - 8 n (n-1)p^*$. Since $\Delta < 0$ when $p^*$ is large for the number of nodes we intend to simulate, we will only consider the cases $p^* = 0.1$ and $p^*=0.5$. 

        When $\Delta > 0$, the root we use is given by 
        \begin{equation}\label{def:mstar}
            m^* = \frac{2n - \sqrt{4n^2 - 8  n (n-1)p^*}}{4}.
        \end{equation}
        For example, when $n= 200$ and $p^* = 0.1$, placing these in the given formula, we find $m^* = 10.501$, which is rounded to 11. Note that for fixed $p^*$, $m^*$ depends linearly on $n$. In particular, when $p^* = 0.5$, $m^* \approx \frac{n}{2}$ and when $p^* = 0.1$, $m^* \approx \frac{n}{20}$.

        4. \textit{Random regular graphs}. Random regular graphs are generated using random\textunderscore regular\textunderscore  graph(d, n) function in NetworkX. Here, $n$ is the number of nodes and $d$ is the common degree. In this case, to achieve an edge density $p^*$, one needs to have $p^* = \frac{\frac{nd}{2}}{\binom{n}{2}}$, which can be rewritten as $p^* = \frac{d}{n-1}$. Thus, given $n$ and $p^*$, we choose $d$ appropriately to achieve the given edge density. 

\subsection{Results for the clustering index} 

Now we are ready to present the simulation results. In each fixed selection of parameters, 120-600 sample graphs of each model were generated for varying node counts, and these were averaged to approximate the value of interest. For each plot, the number of nodes were chosen from the set $\{20, 40, ..., 380\}$, and the resulting averages were linearly interpolated. Note that the rewiring probability in the Watts-Strogatz model does not affect the number of vertices nor the edge density; hence, in each of our simulations, we examined Watts-Strogatz graphs with rewiring probabilities of 0.1, 0.3, 0.5, 0.7 and 0.9. In many of our simulations, Barab\'asi-Albert graph and some of the Watts-Strogatz graphs had growth rates significantly larger than the others; thus, we have created additional plots to better examine the remaining graph models in all our plots. 

The following figure contains our results for $\ci_1$. Note that these plots are for the Barab\'asi-Albert model where the parameter $m = m^*$ is chosen according to \eqref{def:mstar} in order to satisfy the edge density condition. In particular, this is in contrast to the standard Barab\'asi-Albert model where $m$ is a fixed positive integer. One would naturally obtain a distinct growth rate for the standard case. We will have a more detailed relevant discussion below when we focus on the degree index.

\begin{figure}[H]
\begin{center}
    \includegraphics[width=8.1cm]{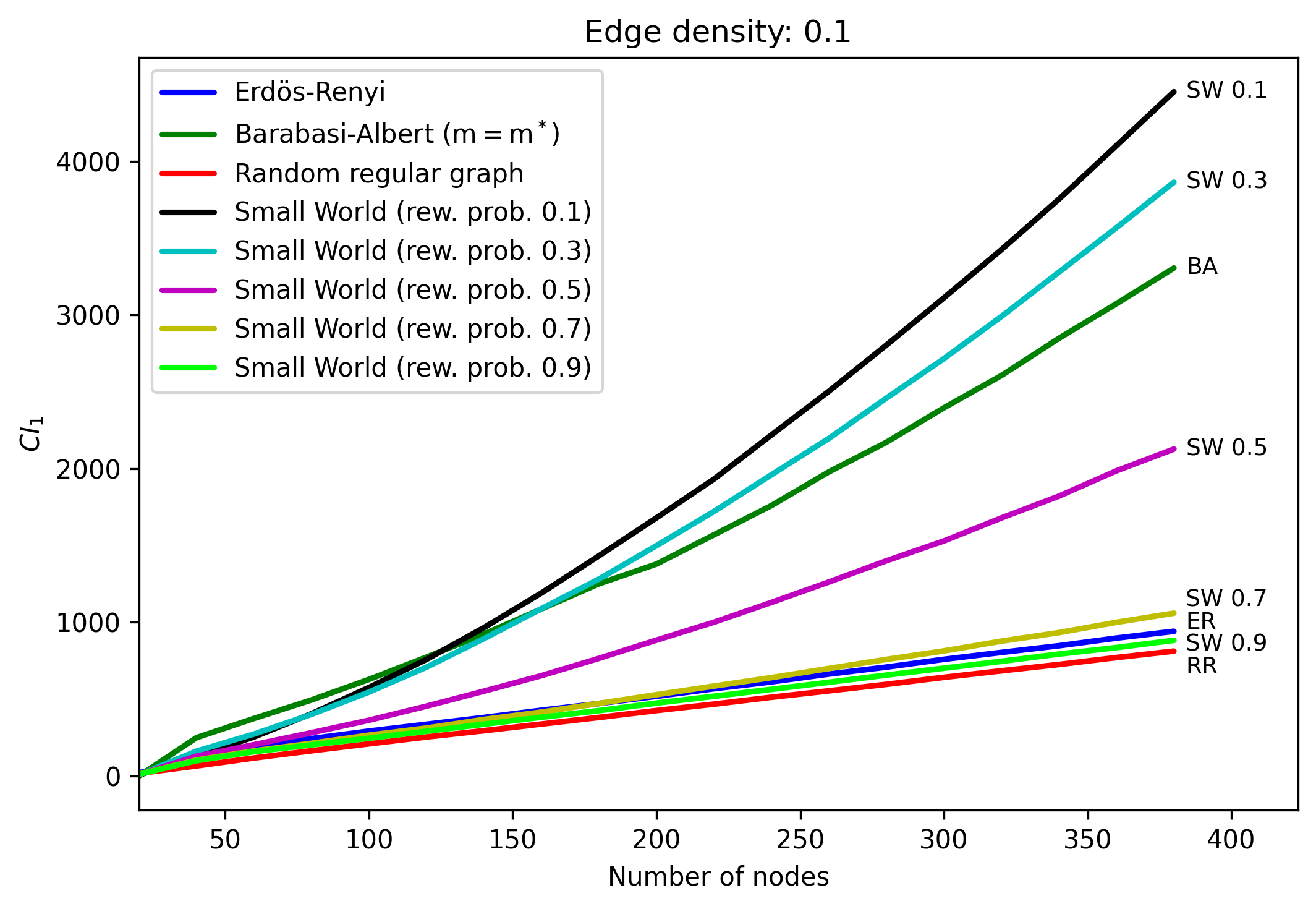} \includegraphics[width=8.1cm]{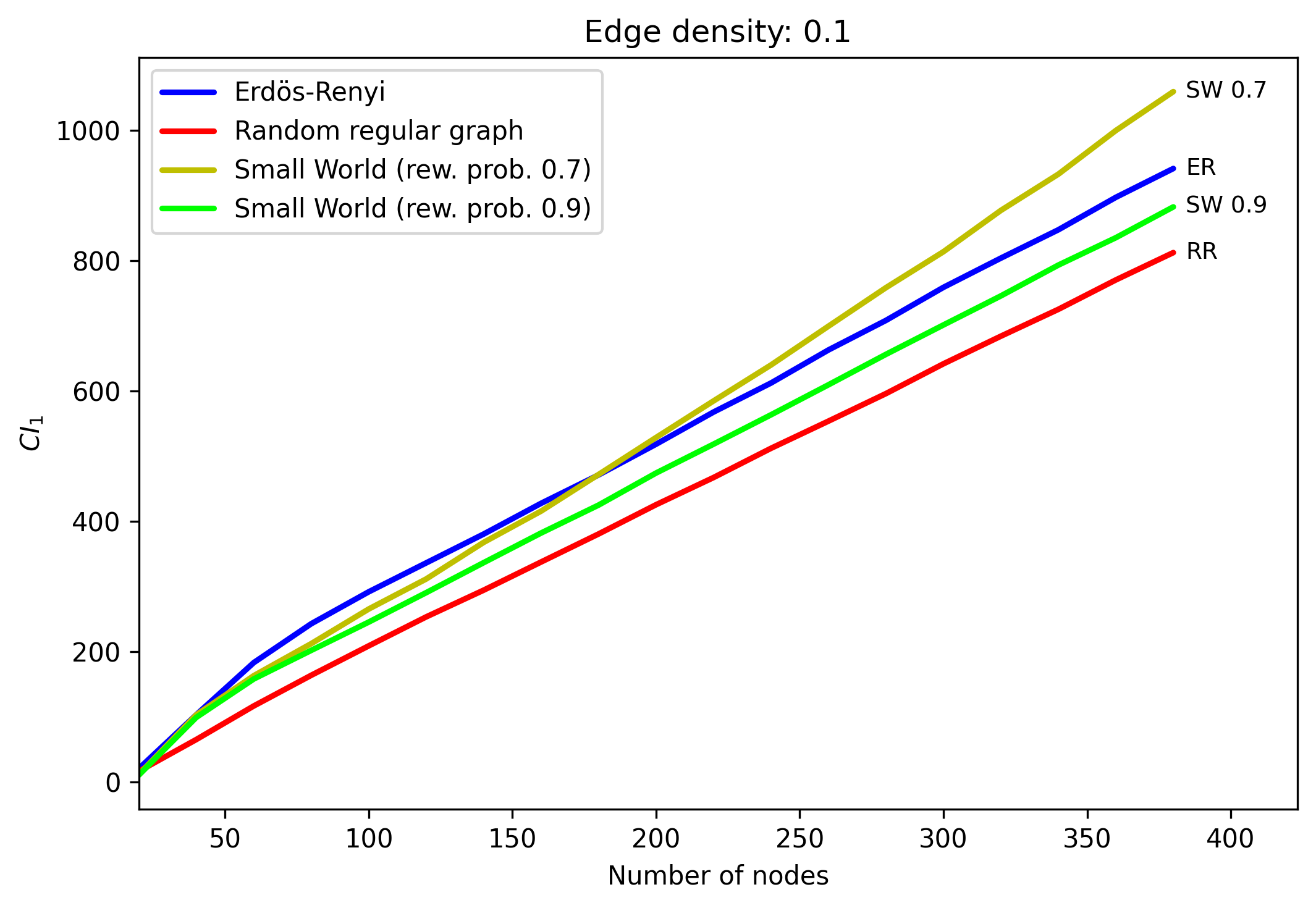}
    \includegraphics[width=8.1cm]{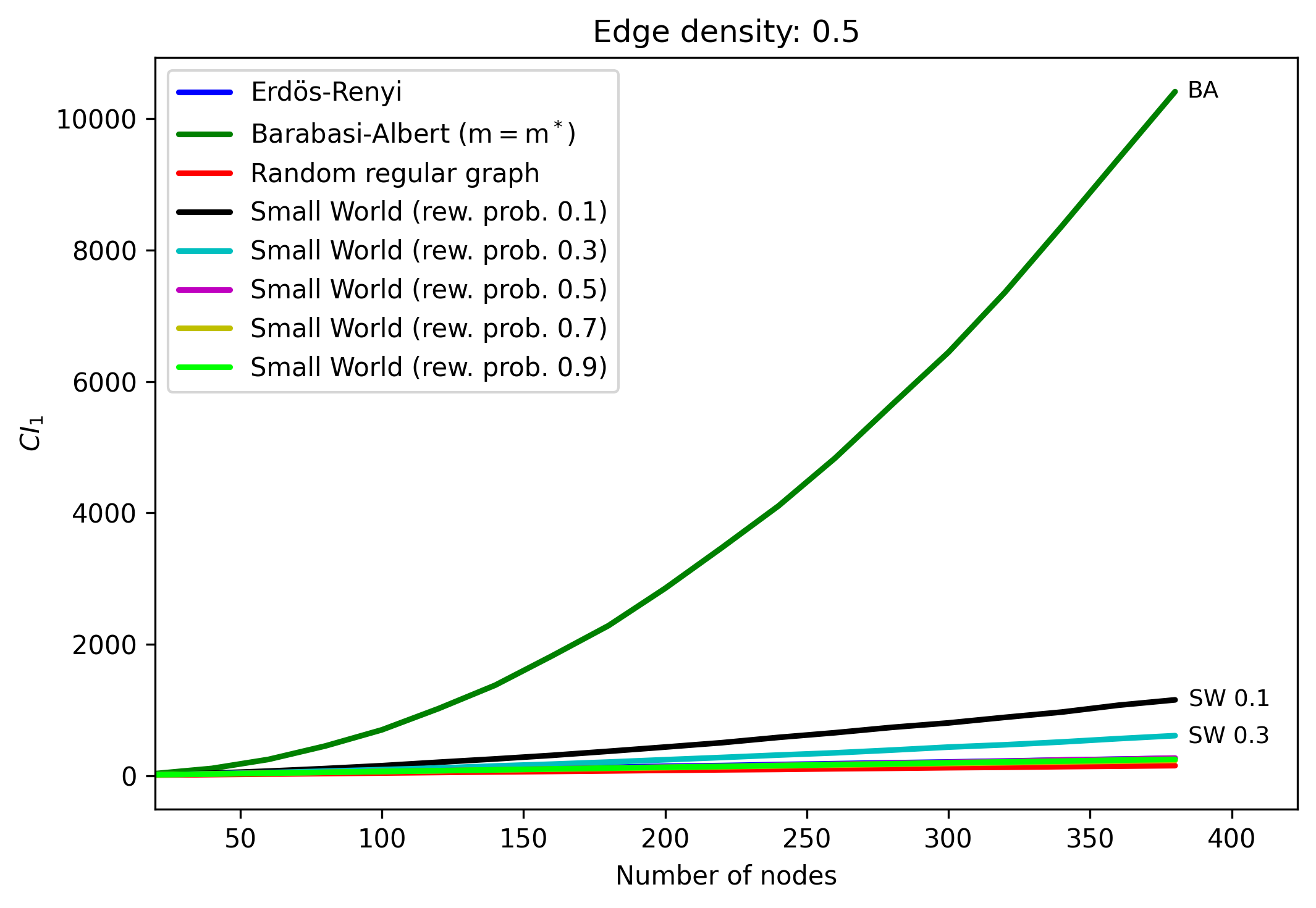} \includegraphics[width=8.1cm]{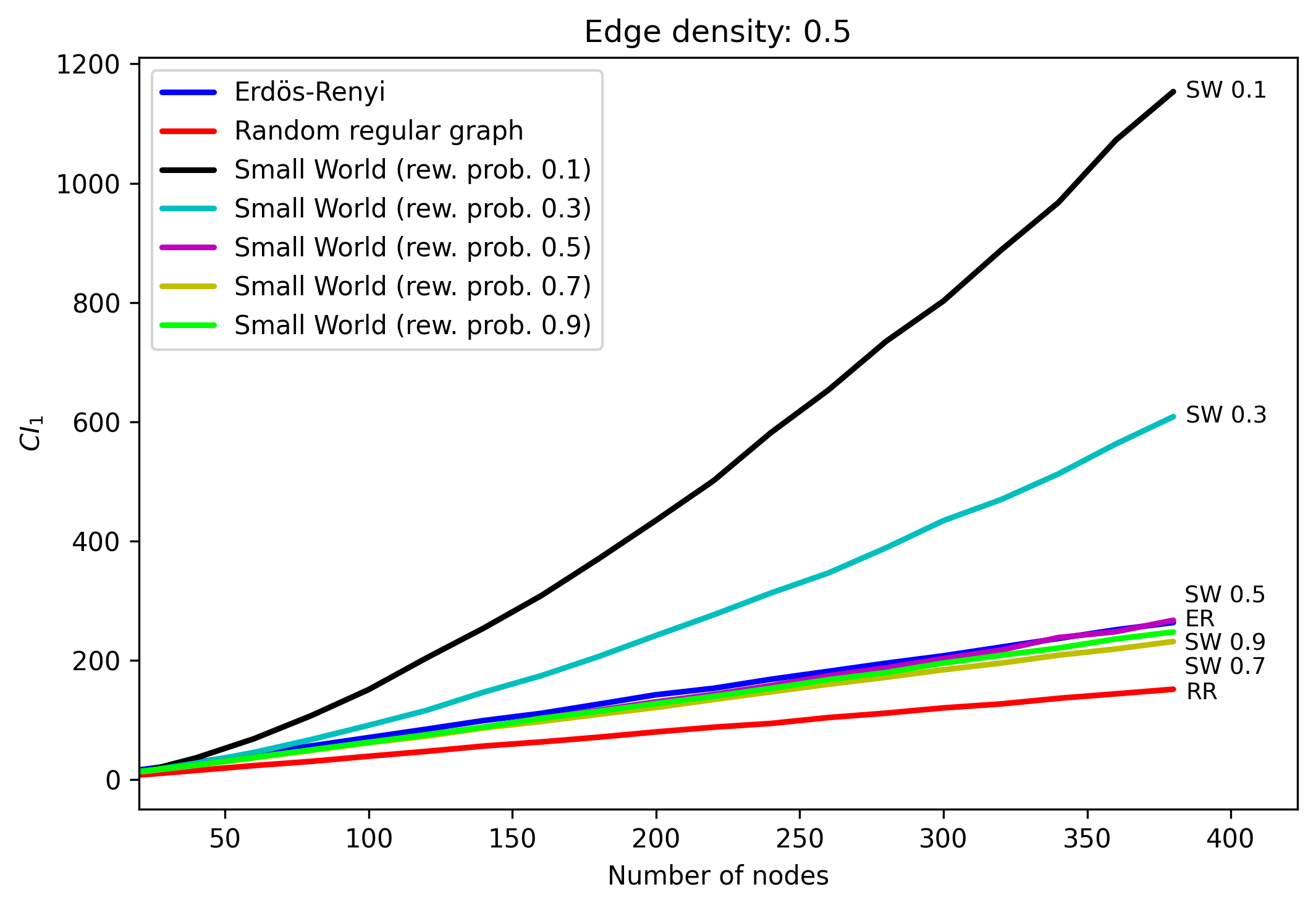}
    \caption{Average $\ci_1$ values of the Erd\H{o}s-R\'enyi graph, Barab\'asi-Albert graph, random regular graph and Watts-Strogatz graphs with edge densities 0.1 and 0.5.}
\end{center}
\end{figure}

\vspace{-0.2in}

Our expectation from our theoretical findings was that the growth of $\ci_1$ for Erd\H{o}s-R\'enyi graphs was sublinear in $n$, which is verified in our simulation. One would also expect that, as the rewiring probability increases, the Watts-Strogatz graphs exhibit increasingly similar behavior to Erd\H{o}s-R\'enyi graphs in terms of clustering, which seems to be the case.


For the Barab\'asi-Albert model with $m=m^*$, when the edge density is 0.5, both $\ci_1$ and $\ci_2$ seem to be of order $n^2$. This is caused by the vastly different clustering behavior of the leaves of the initial star graph and the other nodes. \footnote{Since $m \approx n/2$ in this case, the initial star graph has $m-1 \approx n/2 -1$ leaves. Our simulation results show that, when n is large, the local clustering coefficients of these leaves are significantly higher compared to the clustering coefficients of the other nodes on average. This causes $\mathbb{E}[\ci_1] \geq K_1 \binom{n}{2}$ and $\mathbb{E}[\ci_2] \geq K_2 \binom{n}{2}$ for positive constants $K_1, K_2$.}



\begin{figure}[H]
\begin{center}
    \includegraphics[width=7.0cm]{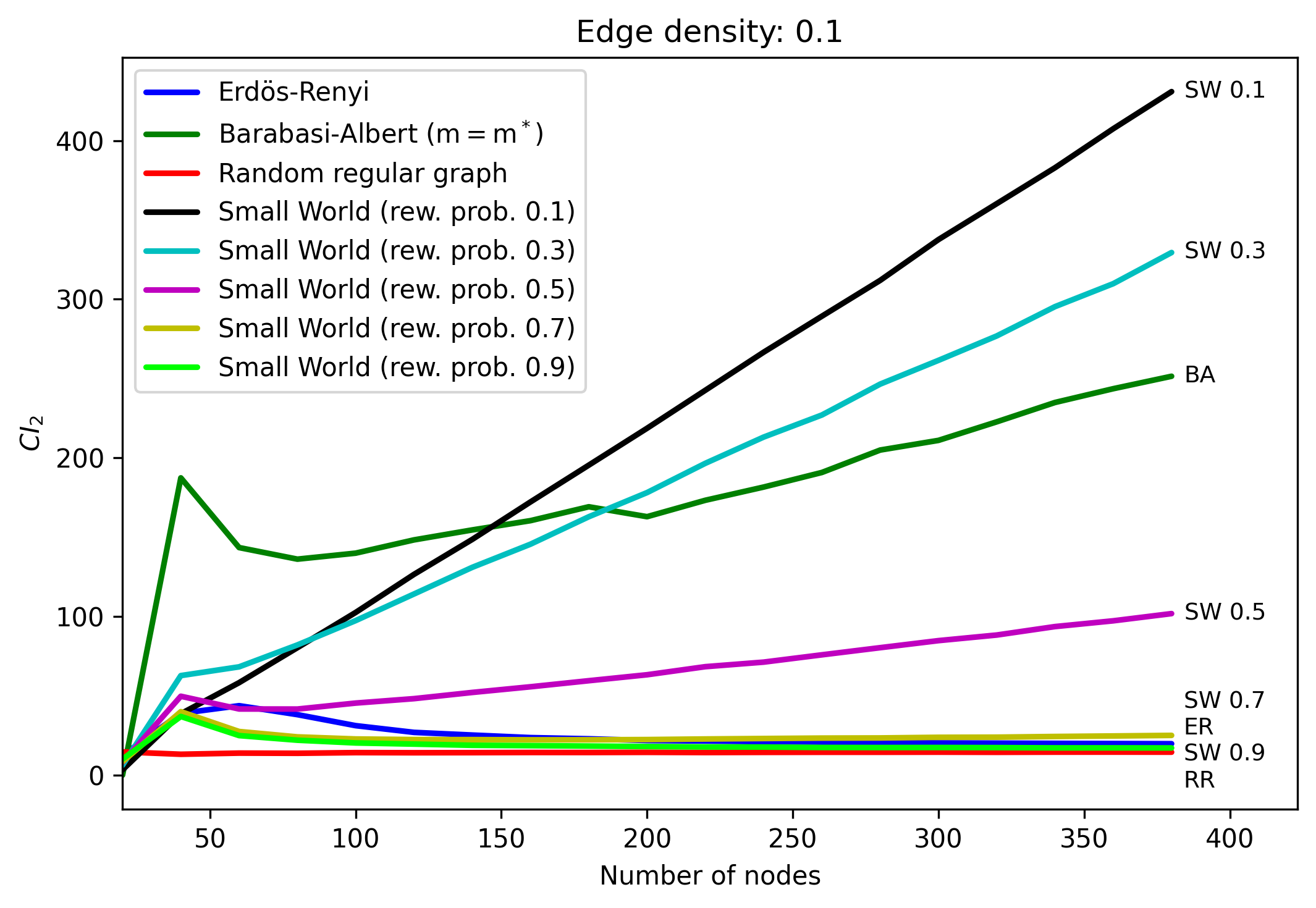} \includegraphics[width=7.0cm]{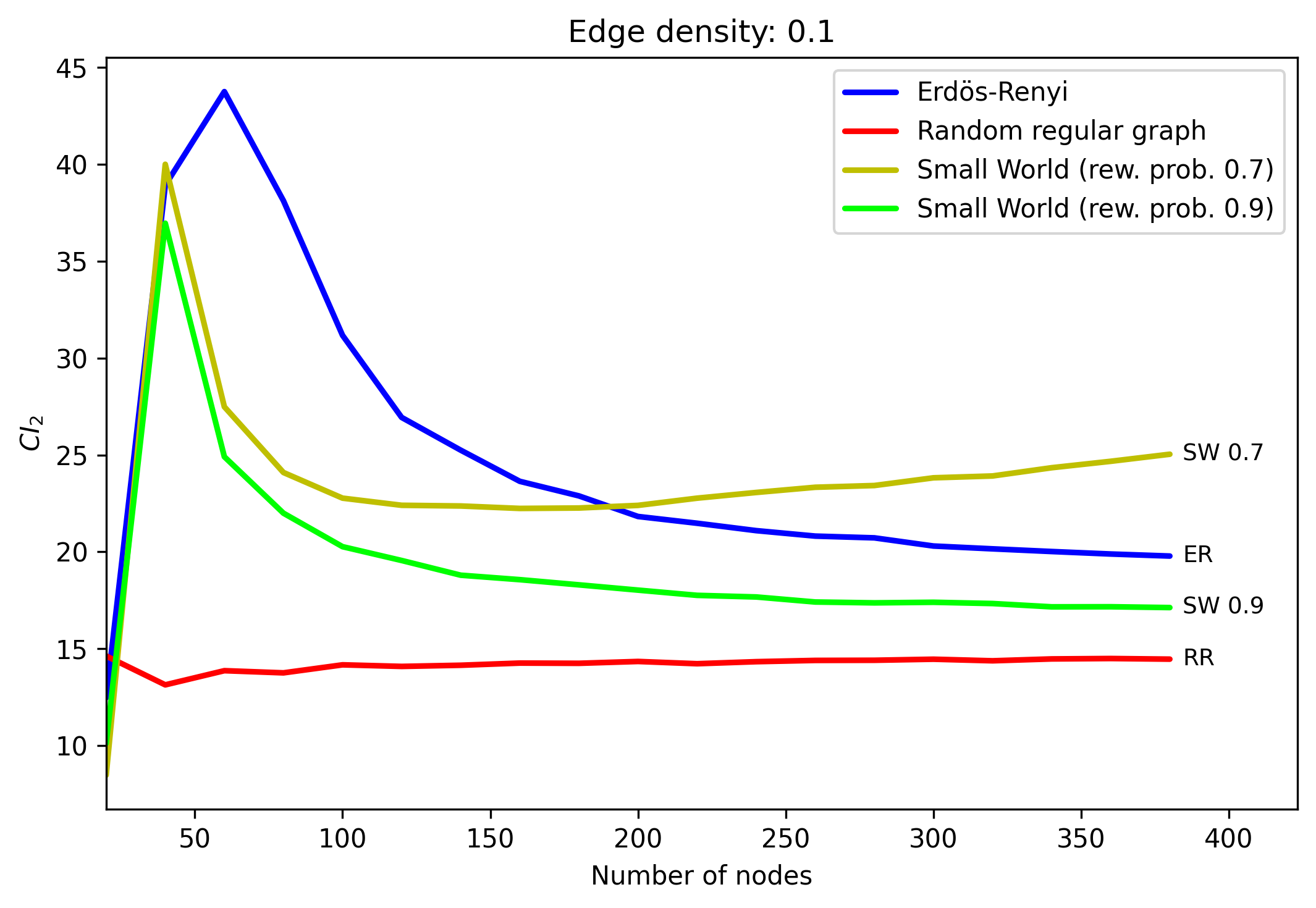}
    \includegraphics[width=7.0cm]{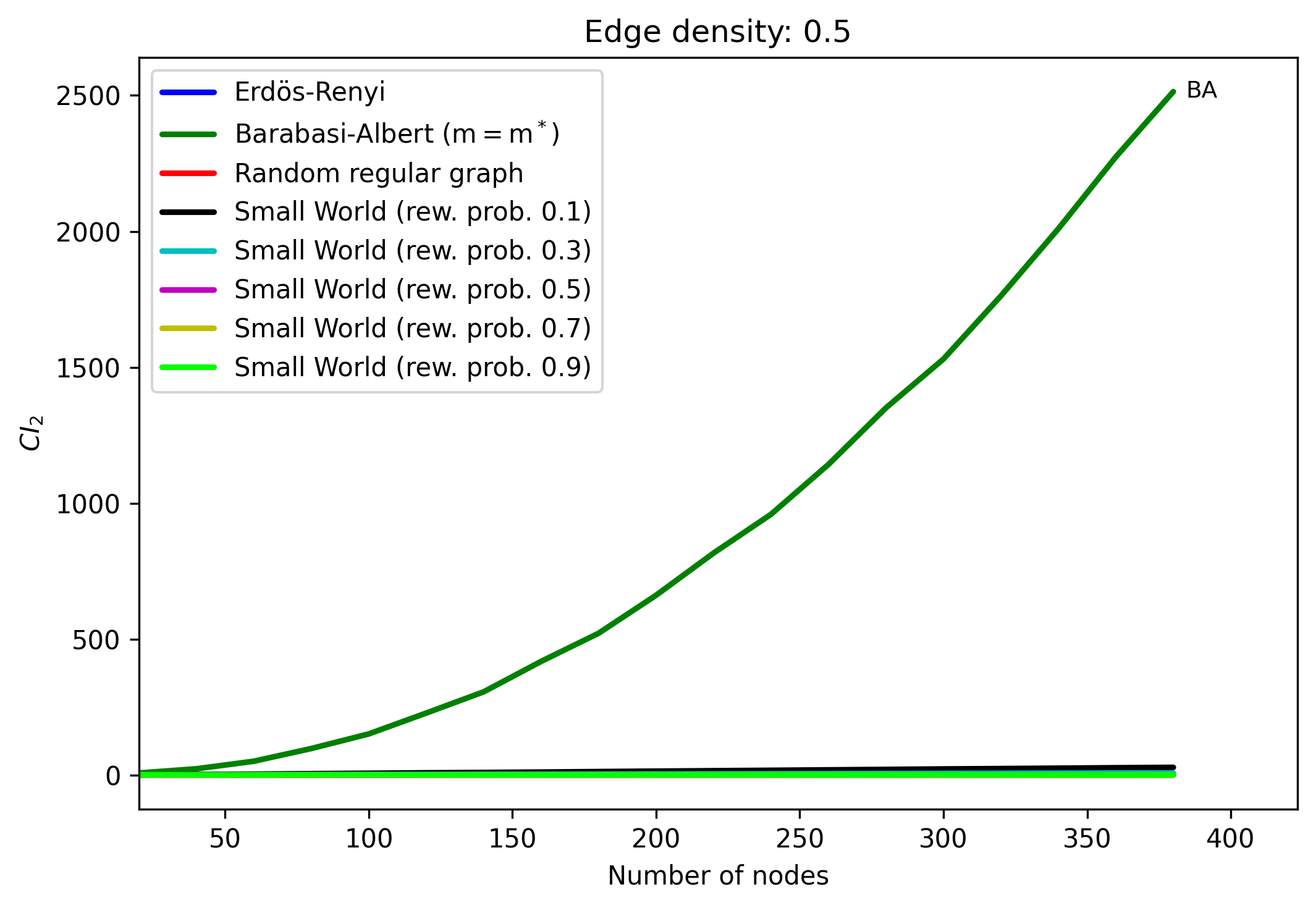} \includegraphics[width=7.0cm]{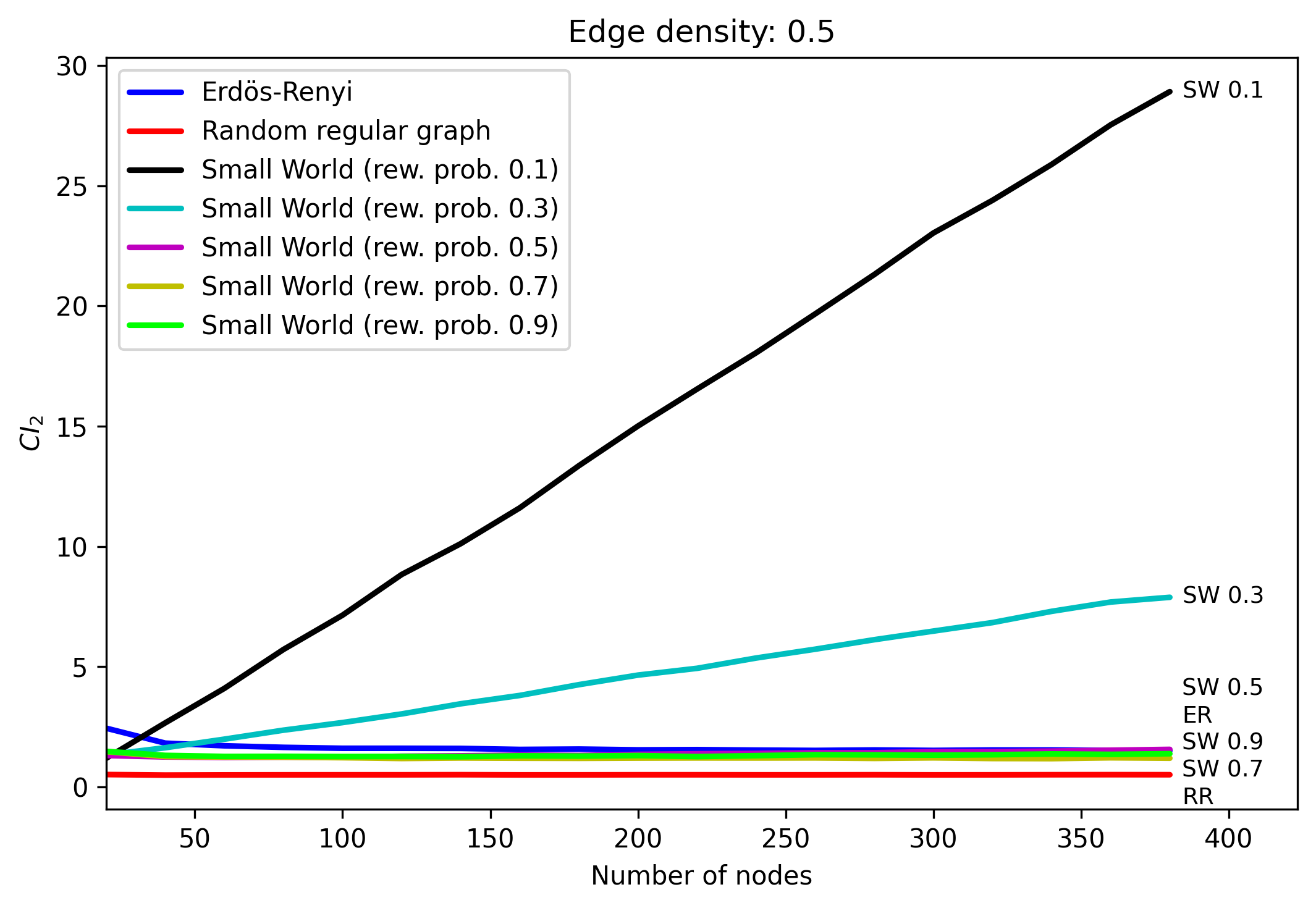}
    \includegraphics[width=7.0cm]{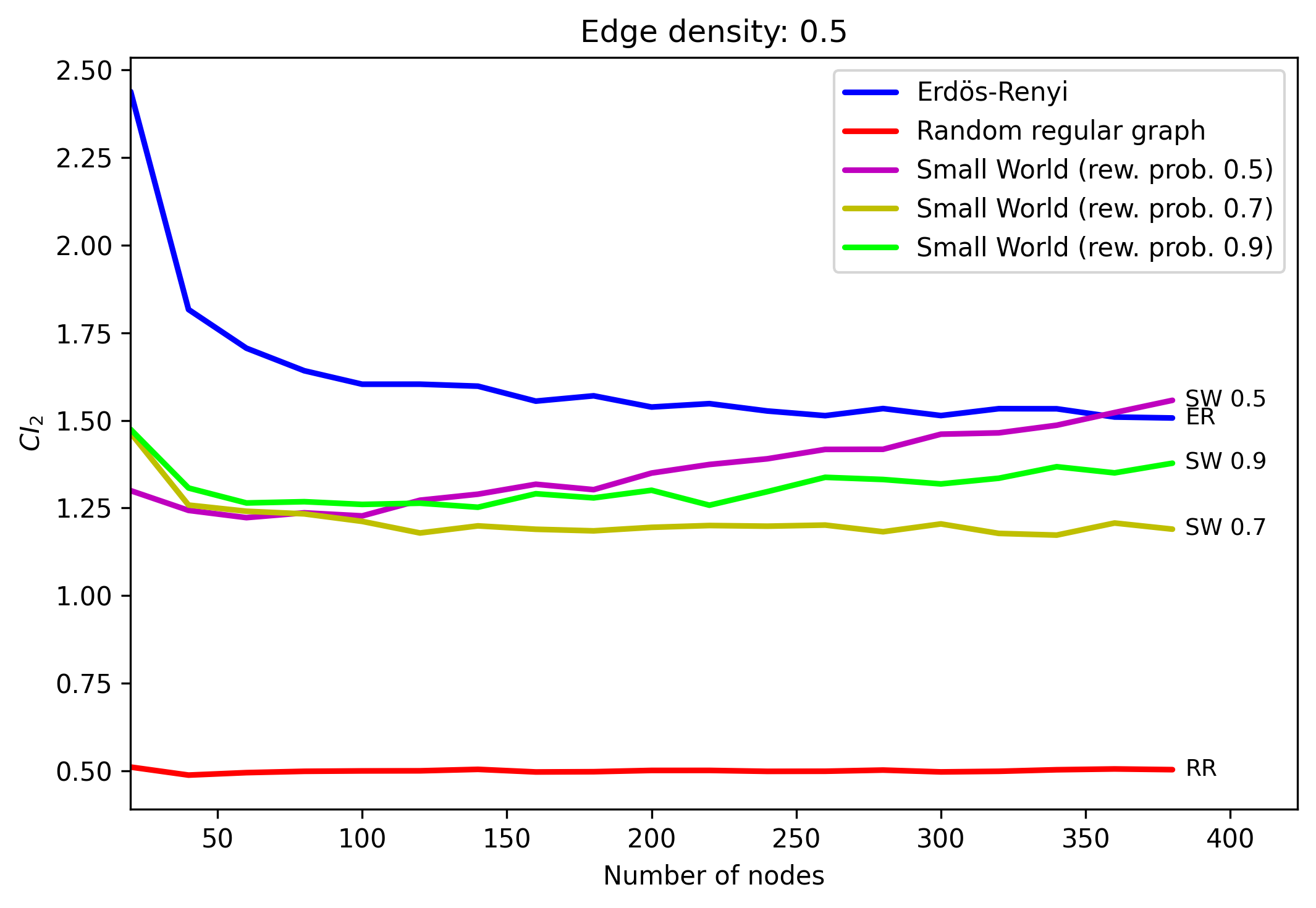}
        \caption{Average $\ci_2$ values of the Erd\H{o}s-R\'enyi graph, Barab\'asi-Albert graph, random regular graph and Watts-Strogatz graphs with edge densities 0.1 and 0.5.}
\label{fig:ci2}
\end{center}
\end{figure}

Next, in Figure \ref{fig:ci2}, we include the comparison of $\ci_2$ for the four models of interest. In the case of Erd\H{o}s-R\'enyi, we showed that this index is upper bounded by a constant independent of the number of nodes and we conjecture that it is indeed also lower bounded by a constant. This is indeed the case according to our simulations, where $\ci_2$ can be observed to approach $\frac{2(1-p)(1-p^2)}{p}$ as $n$ increases. Watts-Strogatz graphs with low rewiring probability and the Barab\'asi-Albert model again dominate the others in this case. Also, as in the case of $\ci_1$, we observe that as the rewiring probability increases, $\ci_2$ of the Watts-Strogatz model approaches that of Erd\H{o}s-R\'enyi model.

\vspace{-0.35in}

It is also interesting in Figure \ref{fig:ci2} to note that many of our random graphs seem to exhibit a peaking behavior around 40-60 nodes when the edge density is 0.1. 
We think that the appearance of these peaks is a result of an equilibrium of two opposing effects:
\begin{itemize}
    \item As $n$ increases, the number of terms appearing in the sum $\ci_2$ also increases.
    \item As $n$ increases, the variance of the distribution of local clustering coefficients decreases.
\end{itemize}
When the edge density is 0.5, the Barab\'asi-Albert model with $m = m^*$ again dominates the others; it seems to exhibit $\Theta(n^2)$ growth while the others are sublinear. Watts-Strogatz model with lower rewiring probabilities again follows Barab\'asi-Albert case in the ranking.

\subsection{Results for the degree index}

\begin{figure}[H]
\begin{center}
    \includegraphics[width=8.0cm]{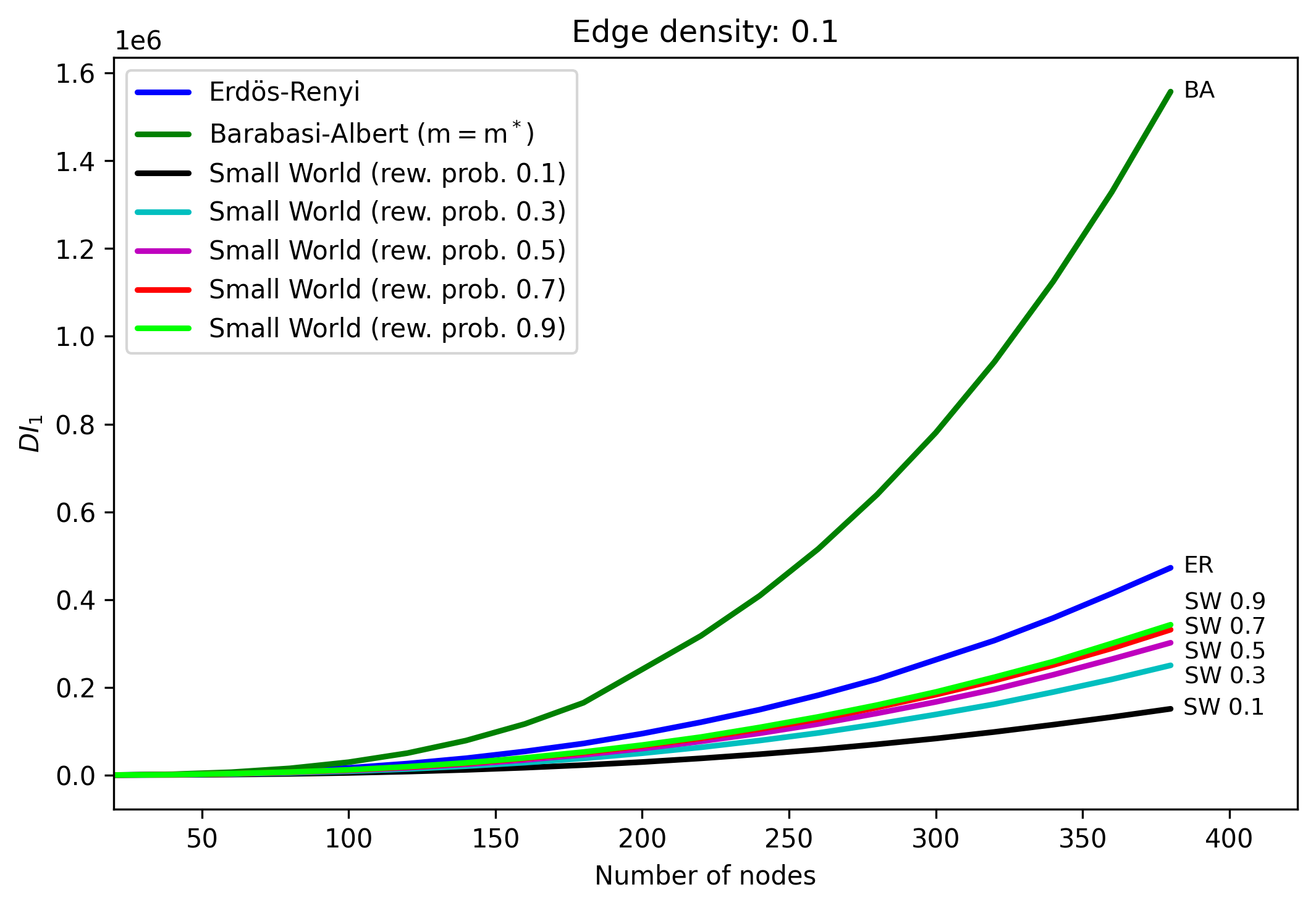}
    \includegraphics[width=8.3cm]{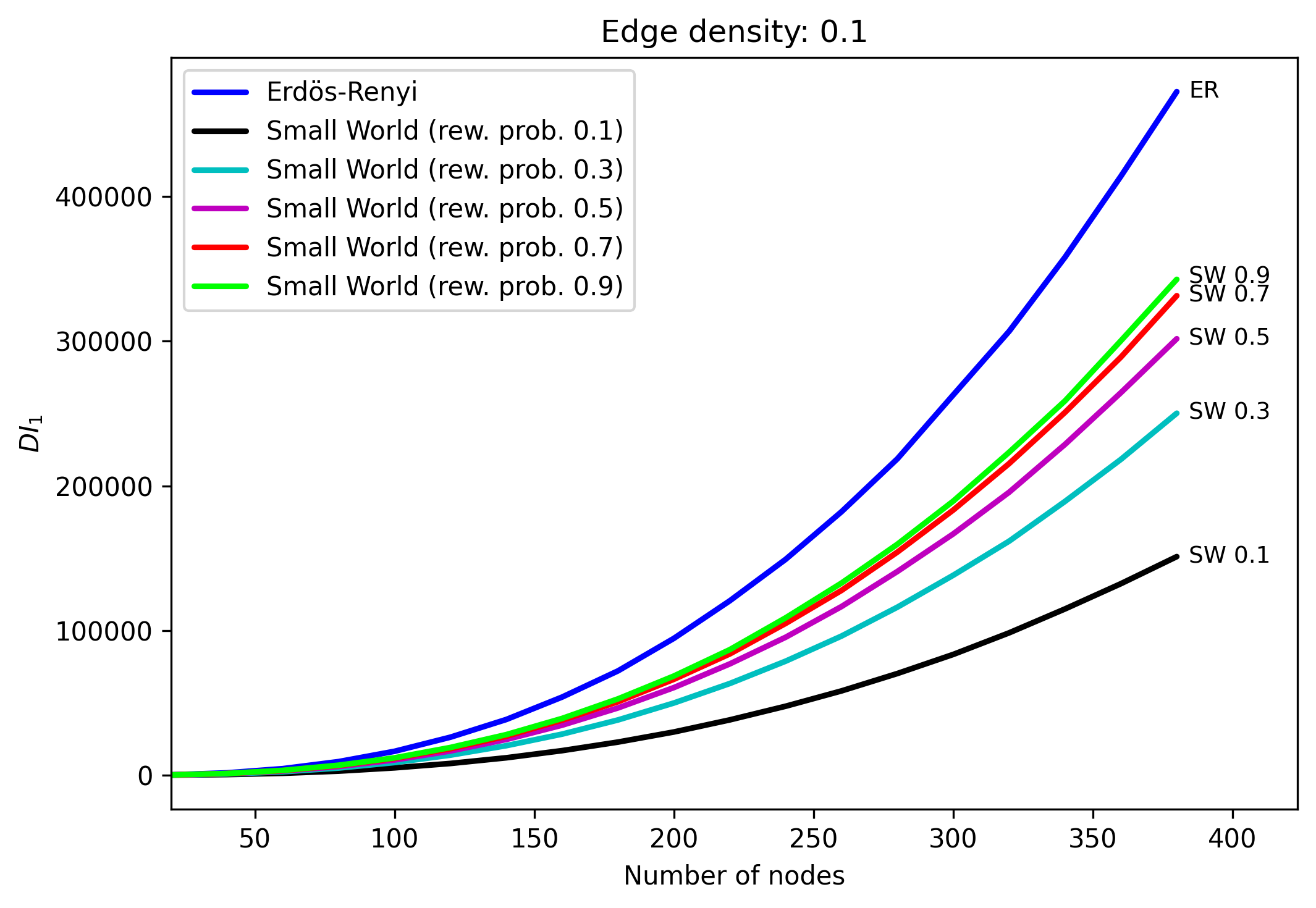}
    \includegraphics[width=8.0cm]{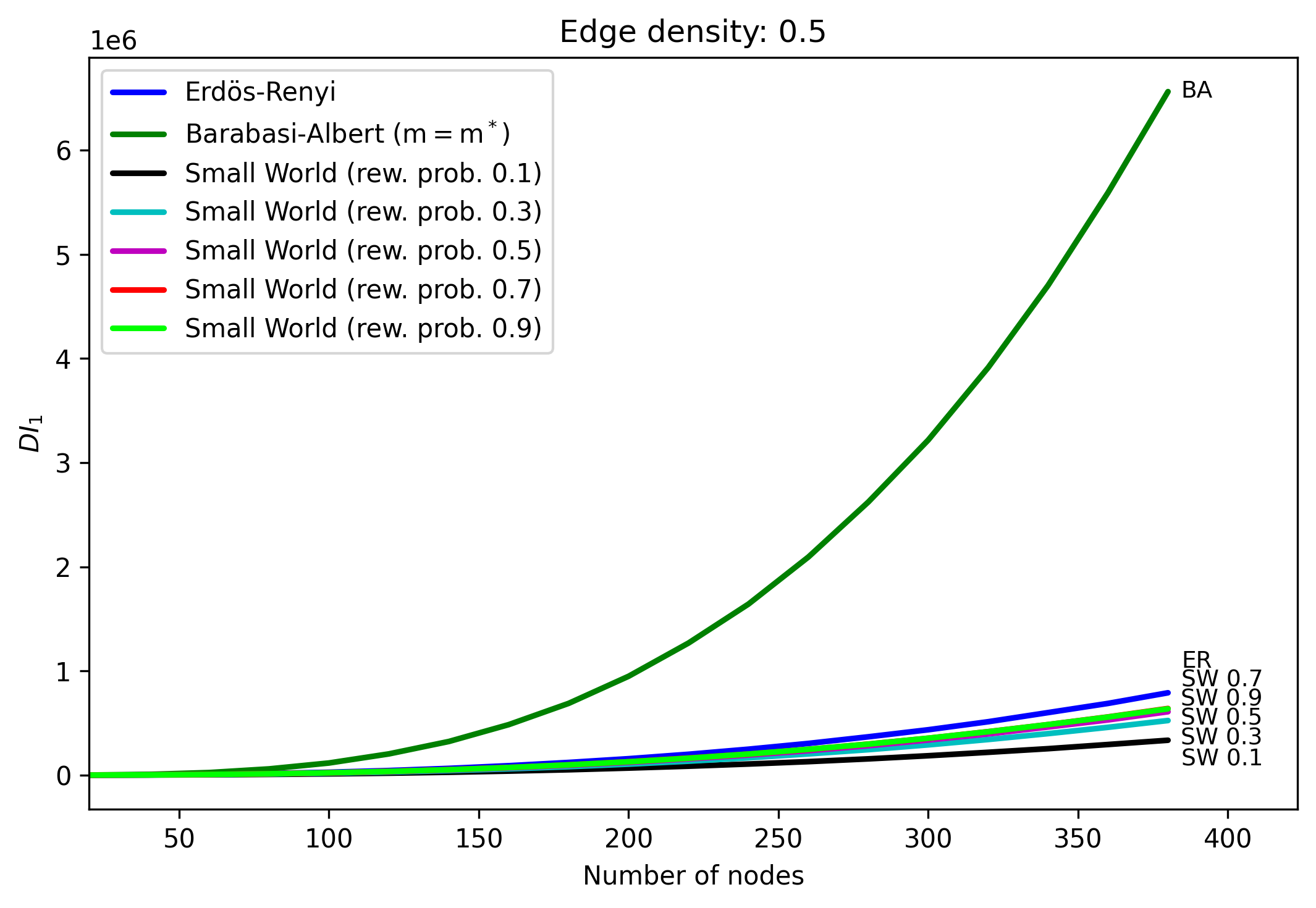}
    \includegraphics[width=8.3cm]{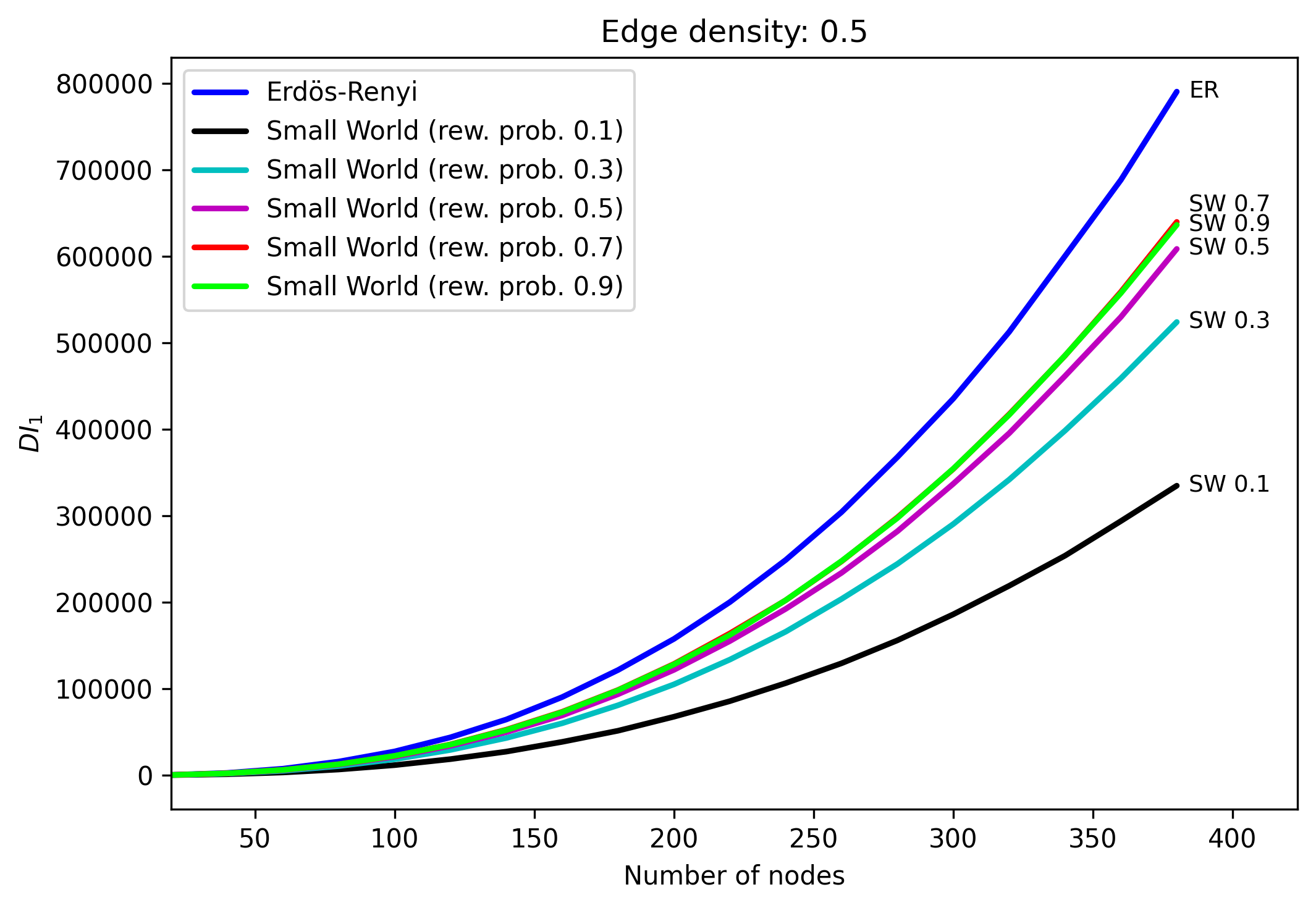}
    \caption{Average $\di_1$ values of the Erd\H{o}s-R\'enyi graph, Barab\'asi-Albert graph, random regular graph and Watts-Strogatz graphs with edge densities 0.1 and 0.5.}
\label{fig:di11}
\end{center}
\end{figure}

In this subsection, we focus on the degree index and compare $\di_1$ and  $\di_2$ for the same four random graph models again for varying values of the number of nodes. Figure \ref{fig:di11} is for $\di_1$ and the edge density in $\{0.1, 0.5\}$. Note that again there are two versions corresponding to each edge density due to the domination of Barab\'asi-Albert model with $m = m^*$. In particular, the ones on the right do not contain the case for the Barab\'asi-Albert model.

\vspace{-0.25in}



We see that as in the case of the clustering index, the small world model behaves similarly to Erd\H{o}s-R\'enyi as the rewiring probability increases. In order to examine the behaviors for different models in more detail, we also include the $\log-\log$ plots for the graphs on the left in Figure \ref{fig:di11}. 

\vspace{-0.13in}
\begin{figure}[H]
\begin{center}
    \includegraphics[width=7.6cm]{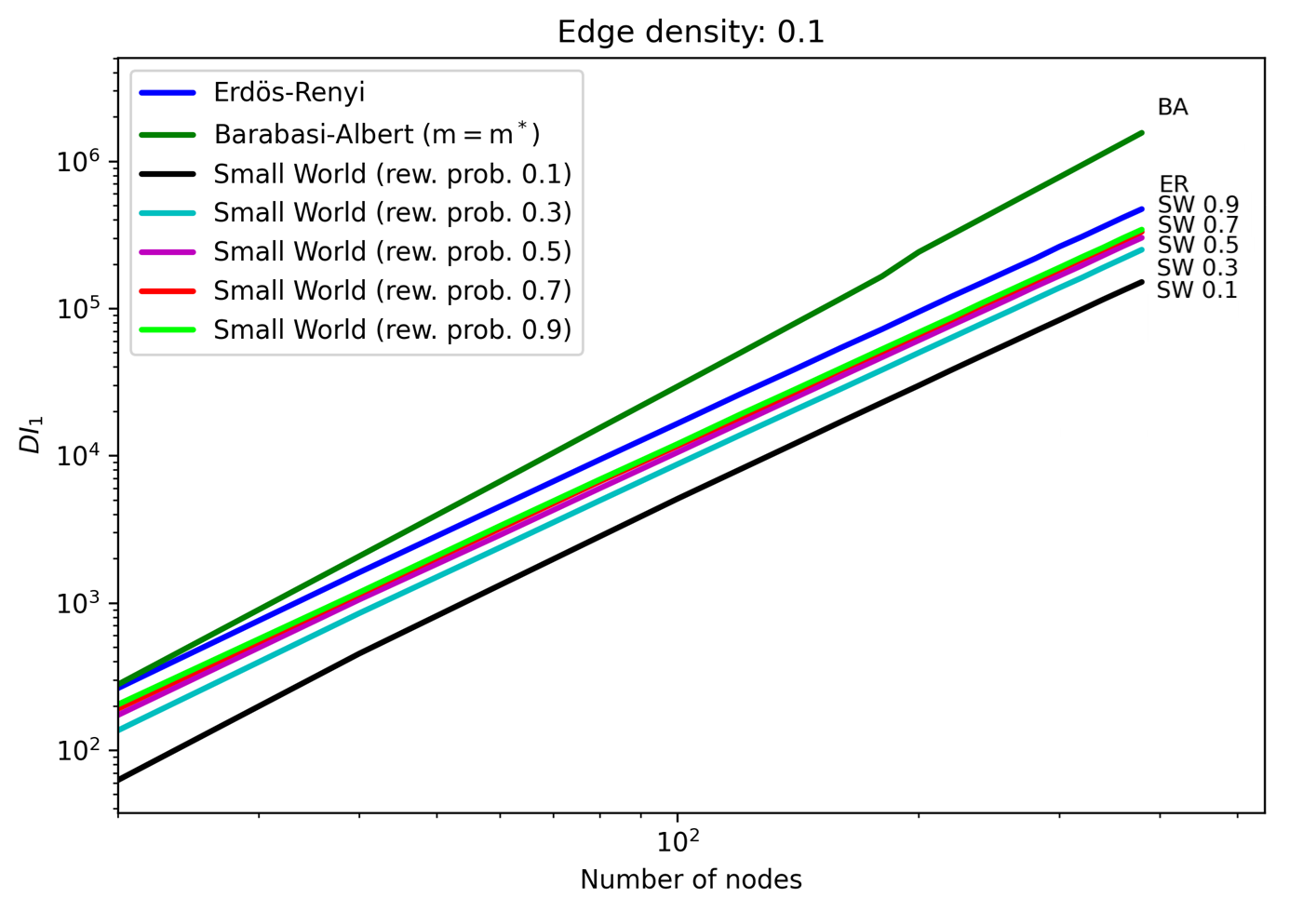}
    \includegraphics[width=7.6cm]{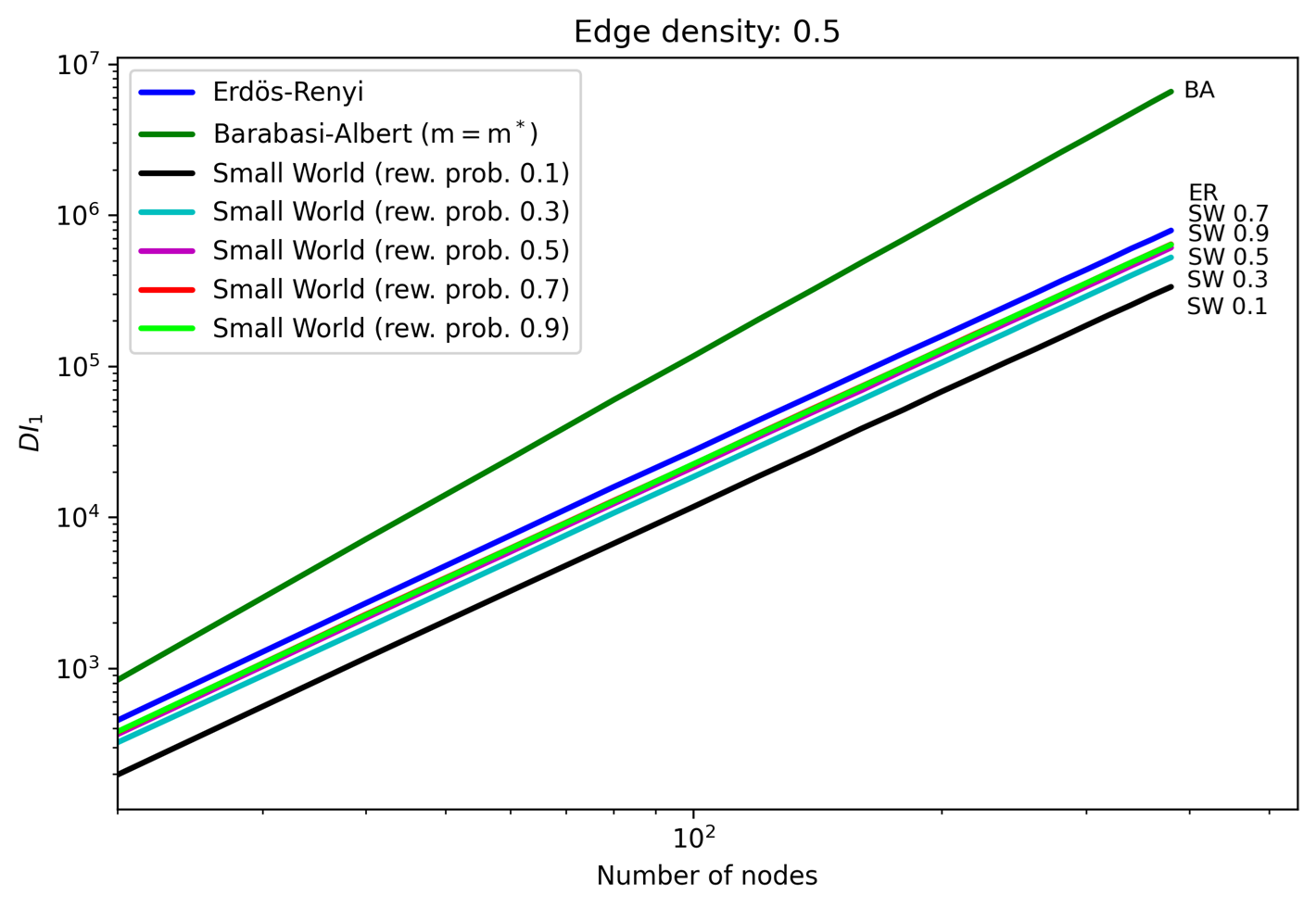}
    \caption{$\log-\log$ plots for the $\di_1$ values of the Erd\H{o}s-R\'enyi graph, modified Barab\'asi-Albert graph, and Watts-Strogatz graphs with edge densities 0.1 and 0.5.}
\label{fig:di1log}
\end{center}
\end{figure}

\vspace{-0.2in}

\begin{figure}[H]
\begin{center}
    \includegraphics[width=7.3cm]{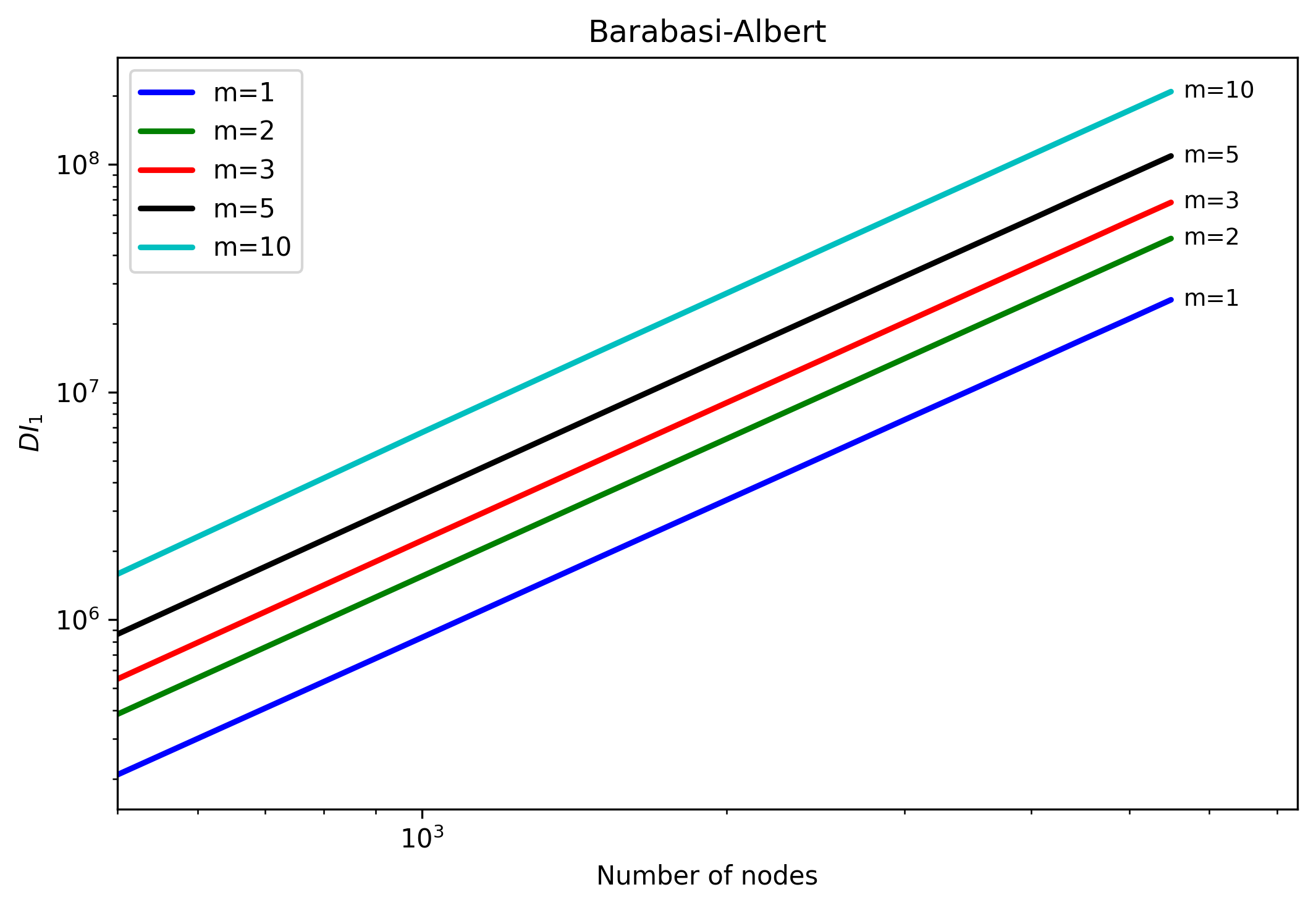}
    \includegraphics[width=7.3cm]{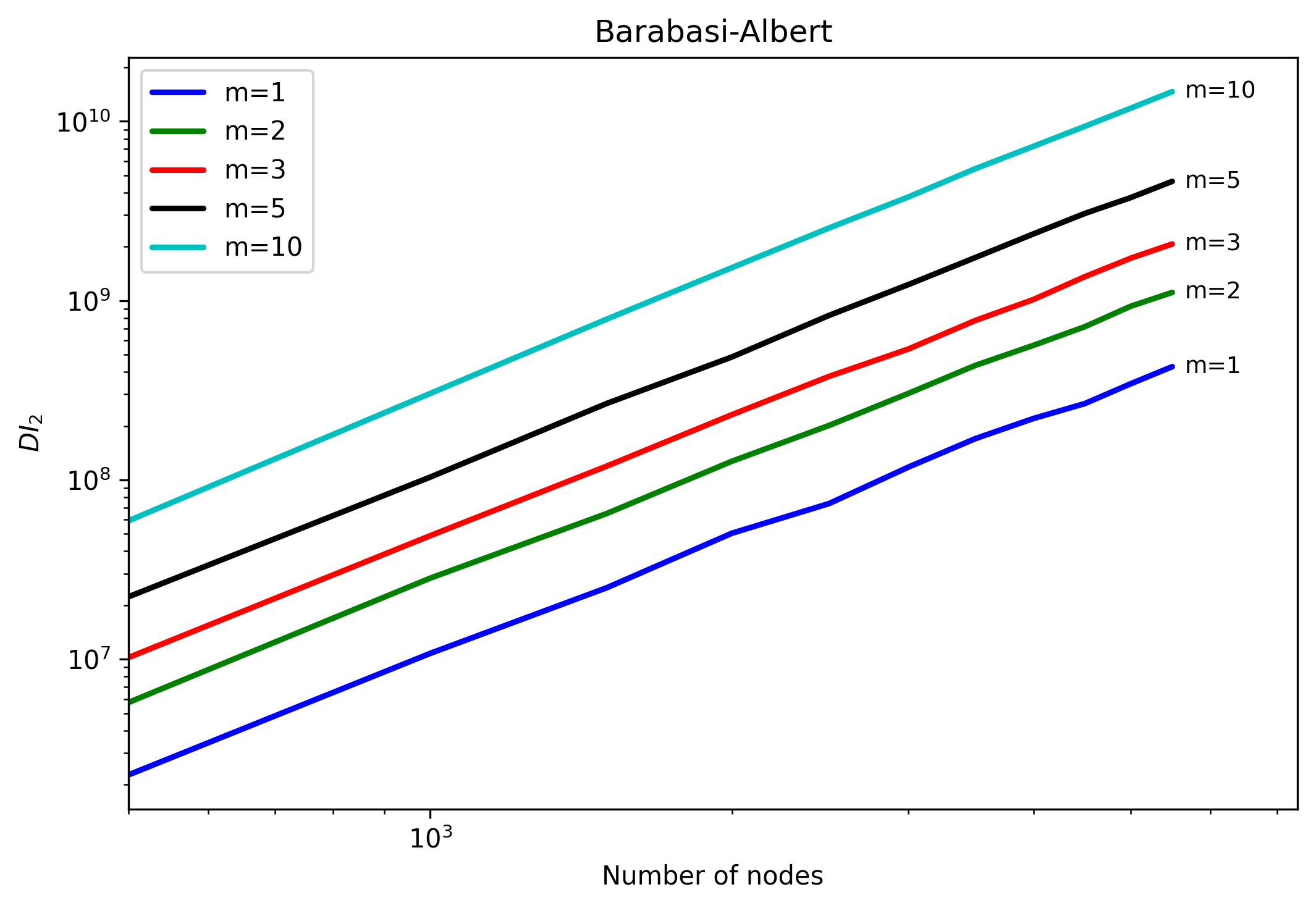}
    \caption{$\log-\log$ plots for the $\di_1$ and $\di_2$ values of the Barab\'asi-Albert graphs with $m \in \{1,2,3,5,10\}$, where $m$ denotes the number of edges to attach from a new node to existing nodes. The number of nodes is from 500 to 6000, with increments of 500.}
\label{fig:dibarlog}
\end{center}
\end{figure}

The $\log-\log$ plot clearly indicates that the growth of $\mathbb{E}[\di_1]$ for the preferential attachment model we consider is of higher order than the growth of other cases. The simulations we have and the heuristic arguments we use suggest that the growth rate in this case is of order $\Theta(n^3)$, which is strictly larger than the one for Erd\H{o}s-R\'enyi graphs (which is $\Theta(n^{5/2})$). However, note that the Barab\'asi-Albert model we consider is different from the standard Barab\'asi-Albert model in the following sense. Here, in order to have a fixed edge density for comparison among distinct models, the corresponding parameter $m$ increases as a function of $n$ as explained in Section \ref{sec:otherrgm}; see \eqref{def:mstar}. Thus, if we consider the standard Barab\'asi-Albert model where $m$ is fixed, then we would obtain a different behavior. Indeed, 
we were informed by an anonymous reviewer that the maximal degree in the standard Barab\'asi-Albert model is of order $\sqrt{n}$; see, for example, page 280 of \cite{hofstadt:24}. From this, one immediately obtains a trivial upper bound of order $n^{5/2}$ for $\di_1$ for fixed $m$. To have a better understanding, we also simulated the standard Barab\'asi-Albert model. Figure \ref{fig:dibarlog} contains the corresponding $\log-\log$-plot.

\vspace{-0.13in}

\vspace{-0.2in}

The plots in Figure \ref{fig:dibarlog} and the related slope calculations suggest that $\mathbb{E}[\di_1] = \Theta(n^2)$ when $m$ is fixed. This seems natural since in this case the degrees of most nodes tend to be small and close to each other, which then causes $\di_1$ to be small. 
On the other hand, in the Barab\'asi-Albert model with $m = m^*$, \eqref{def:mstar} gives $m^* \sim C n$ for some $C \in (0,1)$. In turn, this causes a significant number of nodes to have degrees comparable to $n$ and a remarkable number of nodes to have significantly smaller degrees, 
with the result $\mathbb{E}[|d_i - d_j|]$ being of order $n$, and $\mathbb{E}[\di_1]$ of order $n^3$. 

\begin{figure}[H]
\begin{center}
    \includegraphics[width=8.1cm]{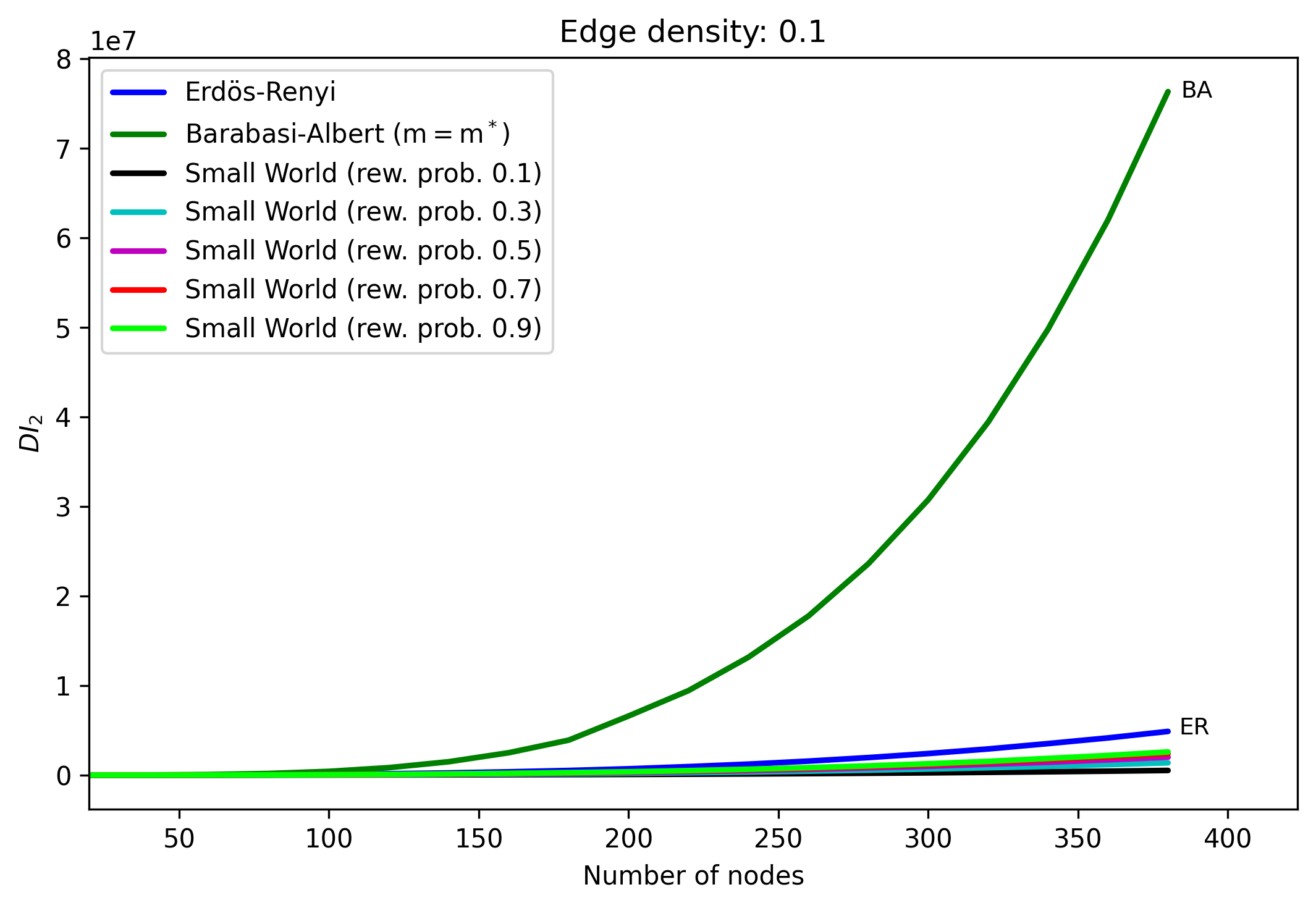}
    \includegraphics[width=8.1cm]{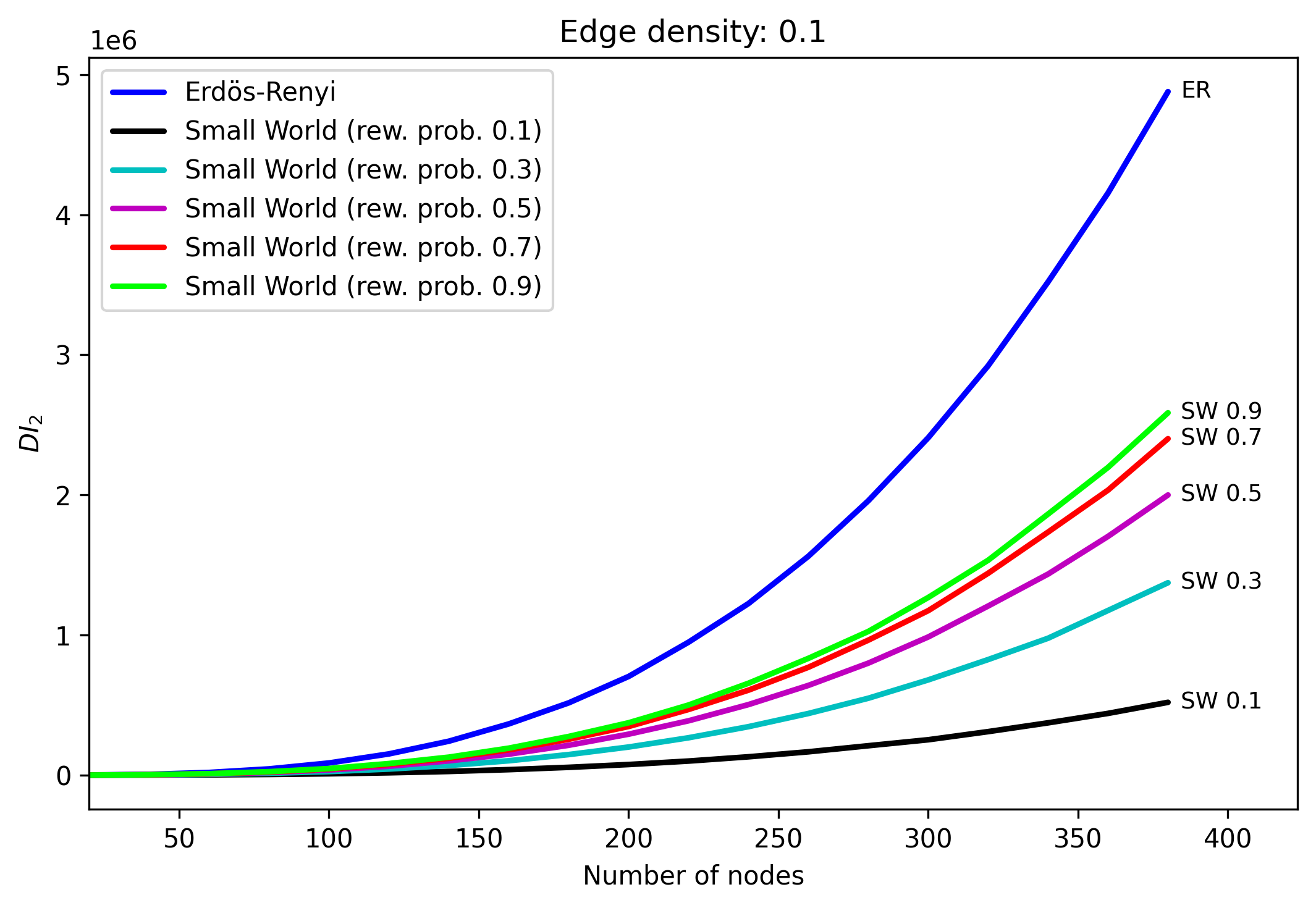}
    \includegraphics[width=8.1cm]{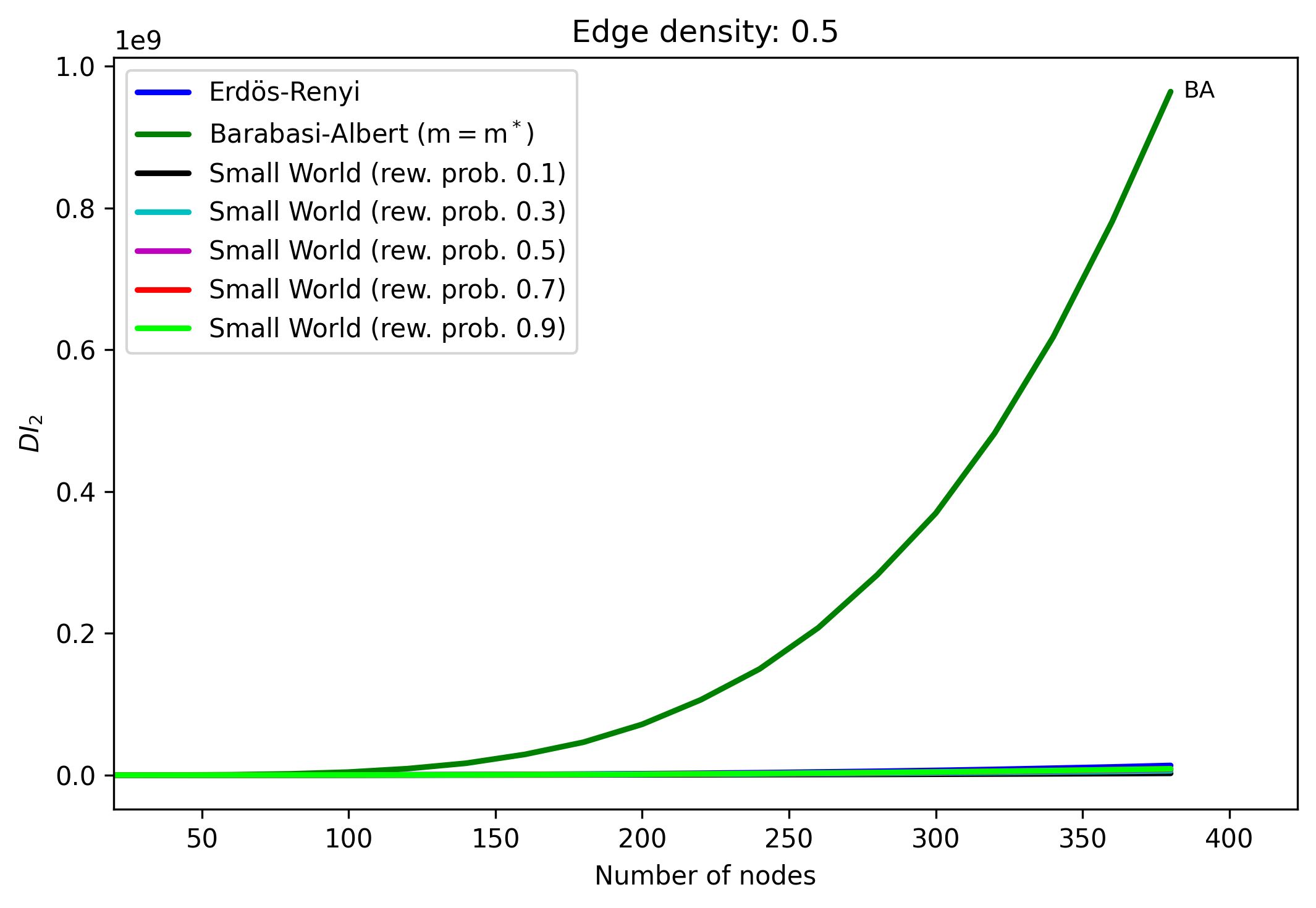}
    \includegraphics[width=8.1cm]{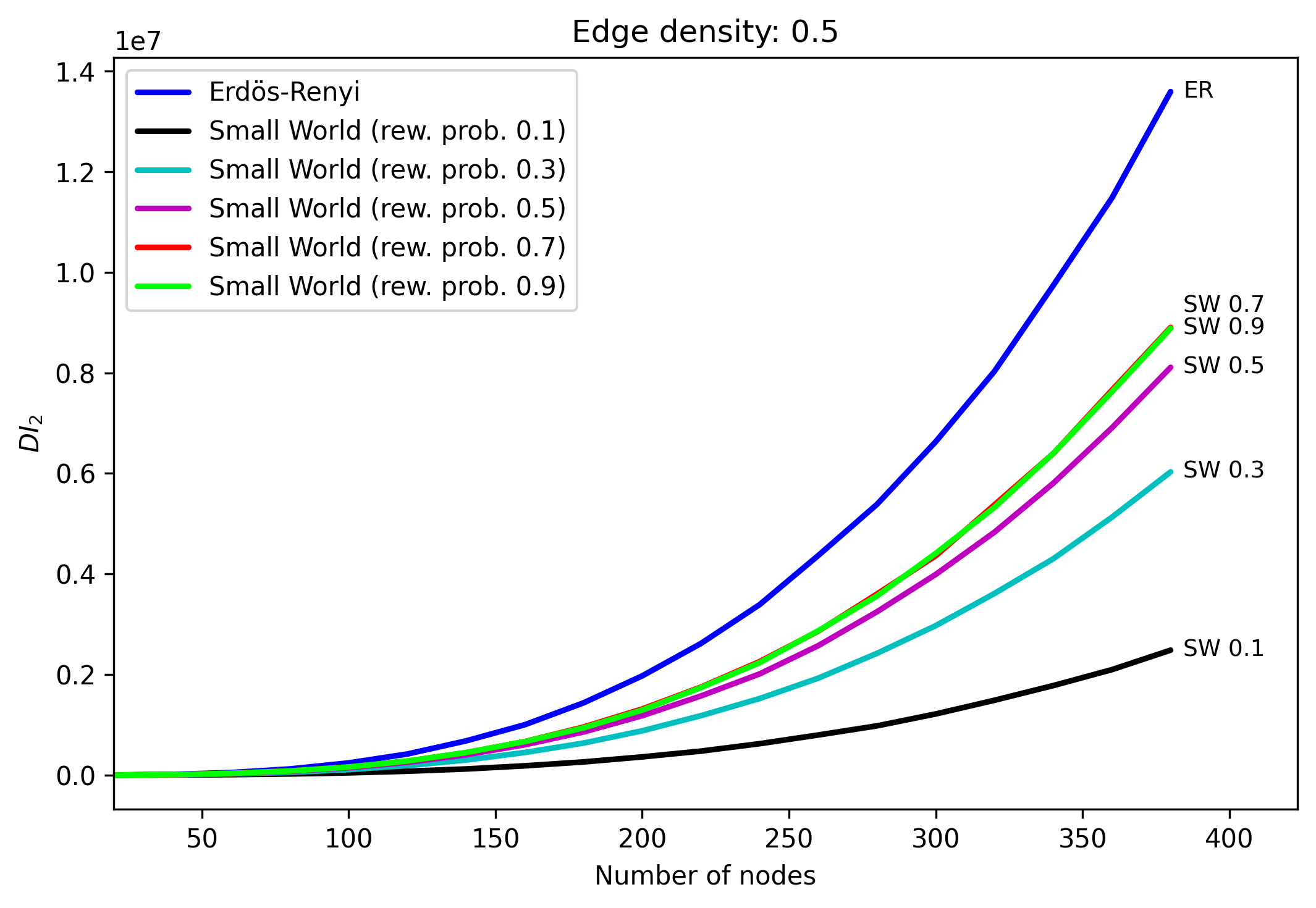}
    \caption{Average $\di_2$ values of the Erd\H{o}s-R\'enyi graph, modified Barab\'asi-Albert graph, random regular graph and Watts-Strogatz graphs with edge densities 0.1 and 0.5.}
\label{fig:di12}
\end{center}
\end{figure}

For the $\di_2$ case, the simulation results are summarized in Figure \ref{fig:di12}. 
The observations here are similar to the $\di_1$ case, and we do not detail these to avoid repetition.


\vspace{-0.2in}

We close this section by noting that we also analyzed the $\log-\log$ plot for the case $\di_2$ for the Barab\'asi-Albert model with $m = m^*$, and that the related observations suggest $\mathbb{E}[\di_2] = \Theta(n^4)$. 

        \section{Conclusion}\label{sec:conclusion}

        In this manuscript, in addition to the well-known degree index (or, irregularity), we introduced the concept of clustering index and studied the first moment of both in a random graph framework. Focusing on Erd\H{o}s-R\'enyi graphs and the degree case, we were able to obtain either exact or asymptotic expressions for the expected degree indices of interest. Regarding the clustering index case, we obtained a linear upper bound for $\mathbb{E}[\ci_1(\mathcal{G})]$, and also showed that $\mathbb{E}[\ci_2(\mathcal{G})]$ is bounded above by a constant independent of $n$.  Although we had a theoretical and simulation-based related analysis for random graphs, we have not done a detailed study of the indices of interest for real-life complex networks and we intend to do so in a subsequent work. 

        One main motivation for studying these indices in our case is their possible use in artificial intelligence classification algorithms as features. Characteristics such as average clustering and the average degree are already used in various such problems, and we intend to do an experimental study in future work to see whether the use of degree and clustering indices provides improvements in performance in classification problems. A second related path is in the field of finance, more particularly in the detection of financial crises. In a recent work \cite{hetal:23}, it is argued that a certain graph irregularity measure, namely the Laplace energy, can be used to detect economic crises. Recently, we did some preliminary experiments to see whether the indices introduced in this paper could also be considered for such purposes. 
        The initial analysis indicated that both are useful in this direction, but that the results for degree index are more statistically significant. We plan to conduct a rigorous analysis on whether degree and clustering indices (as certain types of network irregularities) can also be considered as a precursor of financial crises in upcoming work.


        The last direction we are willing to follow is in computing the exact asymptotics of $\ci_1$ and $\ci_2 $ in the case of Erd\H{o}s-R\'enyi graphs. For example, focusing on $\ci_2$,  we were able to obtain a constant upper bound, and simulations suggested that this index indeed behaves like a constant as the number of nodes increases. It would be interesting to verify this convergence theoretically.

\section*{Appendix: Proof of Lemma \ref{propn:bin}}

(i) Let $X_n, Y_n$ be independent binomial random variables with parameters $n$ and $1/2$. 
         We claim that $\mathbb{E}|X_n - Y_n| = \frac{n \binom{2n}{n}}{2^{2n}}$. For this purpose, letting $Z_n$ be an independent binomial random variable with parameters $n$ and $1/2$, observe that  $$|X_n - Y_n| =_d |X_n - n + Z_n | = |X_n +Z_n - n| = |T_n - n|,$$ where $T_n = X_n + Z_n$ is a binomial random variable with parameters $2n$ and $1/2$. Then, $$\mathbb{P} (|X_n - Y_n| = k) = \mathbb{P} (|T_n  -n| = k) = 2 \mathbb{P} (T_n = n +  k) =  \binom{2n}{n + k} \frac{1}{2^{2n - 1}}.$$
        We have, \begin{eqnarray*}
            \mathbb{E}|X_n - Y_n| = \sum_{k = 0}^n k  \binom{2n}{n + k} \frac{1}{2^{2n - 1}} &=& \frac{1}{2^{2n - 1}} \sum_{k= 1}^n (k - n + n) \binom{2n}{n + k} \\
            &=&   \frac{1}{2^{2n - 1}} \sum_{k= 1}^n (n + k) \binom{2n}{n + k} -  \frac{1}{2^{2n - 1}} \sum_{k= 1}^n n  \binom{2n}{n + k} \\
            &=& \frac{2n}{2^{2n - 1}} \sum_{k=1}^n \binom{2n - 1}{n + k  -1} - \frac{n}{2^{2n - 1} } \sum_{k = 1}^n \binom{2n}{n + k} \\
            &=&  \frac{2n}{2^{2n - 1}} \frac{2^{2n - 1}}{2} - \frac{n}{2^{2n - 1}} \left( \frac{2^n - \binom{2n}{n}}{2} \right) = \frac{n \binom{2n}{n}}{2^{2n}}.
        \end{eqnarray*}
        The assertion that $\mathbb{E}|X_n - Y_n|    \sim \sqrt{\frac{n}{\pi}}$, $n \rightarrow \infty$, now follows from Stirling's approximation. 
        
        (ii)  Let $Z_1, Z_2$ be independent standard normal random variables. We denote the Wasserstein distance between the distributions by $d_W$. Now, letting $\mu_n = \mathbb{E}[X_n]$ and $\sigma_n^2 = \var(X_n)$, recall   that $$d_W \left( \frac{X_n - \mu_n}{\sigma_n} , Z_1\right) := \sup \left\{ \mathbb{E} \left| h \left(\frac{X_n - \mu_n}{\sigma_n}  \right) - h(Z_1) \right|: h: \mathbb{R} \rightarrow \mathbb{R}, |h(x) - h(y)| \leq |x - y| \right\} \leq \frac{D}{\sqrt{n}},$$ for some constant $D$ independent of $n$ (See, for example, \cite{ross}). 
        We in particular  have,  
        $$\mathbb{E} \left|  \frac{X_n - \mu_n}{\sigma_n}  -  Z_1  \right| \leq \frac{D}{\sqrt{n}}.$$ 
        With previous observations, and together with triangle inequality, this gives 
        \begin{eqnarray*}
            \mathbb{E}|X_n - Y_n| &=& \sigma_n \left|\frac{X_n - \mu_n}{\sigma_n} - Z_1 + Z_1 - Z_2 + Z_2 - \frac{Y_n - \mu_n}{\sigma_n} \right| \\ 
            &\leq&  2 \sigma_n \mathbb{E} \left|\frac{X_n - \mu_n}{\sigma_n} - Z_1\right| +  \sigma_n \mathbb{E}|Z_1 - Z_2| \\ 
            &\leq& \frac{2 C \sigma_n}{\sqrt{n}} + \sigma_n \mathbb{E}|Z_1 - Z_2| = \mathcal{O}(1) +  \sigma_n \mathbb{E}|Z_1 - Z_2|. 
        \end{eqnarray*}  
        
        Also, using the triangle inequality in an appropriate way,
        \begin{eqnarray*}
            \mathbb{E}|X_n - Y_n | &=& \sigma_n \left| \left( \frac{X_n - \mu_n}{\sigma_n} - Z_1  \right) + (Z_1 - Z_2) + \left(Z_2 - \frac{Y_n - \mu_n}{\sigma_n} \right)\right| \\
            &=& \sigma_n \mathbb{E}|a + b + c| \geq \sigma_n \mathbb{E}|b| - \sigma_n \mathbb{E}|a| -  \sigma_n \mathbb{E}|c| \\
            &=& \sigma_n \mathbb{E}|Z_1 - Z_2| - 2 \sigma_n  \mathbb{E}|\frac{X_n - \mu_n}{\sigma_n} - Z_1| \\
            &\geq&\sigma_n \mathbb{E}|Z_1 - Z_2| -  \frac{2 \sigma_n D}{\sqrt{n}} = \mathcal{O}(1) +  \sigma_n \mathbb{E}|Z_1 - Z_2|.
        \end{eqnarray*}

        So,  $\frac{1}{\sigma_n} \mathbb{E} |X_n - Y_n |$ converges to $\mathbb{E}|Z_1 - Z_2|$. The result now follows since $\mathbb{E}|Z_1 - Z_2| = 2 / \sqrt{\pi}$. \hfill $\square$

\vspace{0.15in}

\textbf{Acknowledgments.} We would  like to thank the anonymous reviewers whose suggestions and corrections improved the paper significantly.  Parts of this paper were completed at the Nesin
Mathematics Village, the authors would like to thank Nesin Mathematics Village for their kind hospitality. The first author is supported by the BAGEP Award of the Science Academy, Turkey.

        \end{document}